\documentclass{article}
\newtheorem{theorem}{Theorem}[section]
\newtheorem{proposition}[theorem]{Proposition}
\newtheorem{lemma}[theorem]{Lemma}

\newtheorem{example}[theorem]{Example}
\newtheorem{definition}[theorem]{Definition}
\usepackage[utf8x]{inputenc}
\usepackage{geometry}
\usepackage{microtype}
\usepackage{pifont}
\usepackage{fancyhdr}
\usepackage{graphicx}
\usepackage{enumerate}
\usepackage{subfig}

\newcommand{\ep}{$\square$}
\newcommand{\proc}[1]{\noindent\textit{#1}}
\newcommand {\cont}[1]{\mathcal{C}\left(#1\right)}
\newcommand{\ceil}[1]{\left\lceil #1 \right\rceil}
\newcommand{\abs}[1]{\left| #1 \right|}
\newcommand{\per}[2]{\mathrm{Per}_{#1}(#2)}
\newcommand{\quotient}[2]{#1/#2}
\newcommand{\set}[2]{\left\{ #1 \, \left| \, #2 \right.\right\}}

\renewcommand{\hat}{\widehat}
\renewcommand{\tilde}{\widetilde}
\renewcommand{\bar}{\overline}
\newcommand{\deff}[3]{#1\colon #2 \rightarrow #3}

\usepackage{bbm}
\usepackage{dsfont}
\usepackage{bbold}
\usepackage{amsmath,amssymb}

\begin{document}
\title{Toeplitz flows and their ordered K-theory}
\author{Siri-Mal\'en H\o ynes, Norwegian University of Science and Technology (NTNU)}

\maketitle
\begin{abstract}
To a Toeplitz flow $(X,T)$ we associate an ordered $K^0$-group, denoted $K^0(X,T)$, which is order isomorphic to the $K^0$-group of the associated (non-commutative) $C^\ast$-crossed product $C(X)\rtimes_T \mathbb{Z}$. However, $K^0(X,T)$ can 
be defined in purely dynamical terms, and it turns out to be a complete invariant for (strong) orbit equivalence. We characterize the $K^0$-groups that arise from Toeplitz flows $(X,T)$ as exactly those simple dimension groups $(G,G^+)$ 
that contain a noncyclic subgroup $H$ of rank one that intersects $G^+$ nontrivially. Furthermore, the Bratteli diagram realization of $(G,G^+)$ can be chosen to have the ERS-property, 
i.e.\ the incidence matrices of the Bratteli diagram have equal row sums. We also prove that for any Choquet simplex $K$ there exists an uncountable family of pairwise non-orbit equivalent Toeplitz flows $(X,T)$ such that the set of $T$-invariant probability measures $M(X,T)$ is affinely 
homeomorphic to $K$, where the entropy $h(T)$ may be prescribed beforehand. Furthermore, the analogous result is true if we substitute strong orbit equivalence for orbit equivalence, but in that case we can actually prescibe both the 
entropy and the maximal equicontinuous factor of $(X,T)$. Finally, we present some interesting concrete examples of dimension groups associated to Toeplitz flows. 
\end{abstract}
\section{Introduction.}

Toeplitz flows have been extensively studied, both as topological and measure-theoretic dynamical systems, since they were first introduced by Jacobs and Keane in 1969 \cite{JK}. (Note: In spite of
the word ``flow'', which would indicate an $\mathbb{R}$-action, Toeplitz flows are dynamical systems with a $\mathbb{Z}$-action generated by a single homeomorphism.) In this paper
we will exclusively study Toeplitz flows as topological dynamical systems. A Toeplitz flow is a special minimal subshift on a finite alphabet $\Lambda$ (in particular, it is a \emph{symbolic system}) 
defined in terms of a so-called \emph{Toeplitz sequence} on $\Lambda$, the latter being associated to certain arithmetic progressions. (Incidentally, the reason the name ``Toeplitz'' has been 
attached to these systems is that Toeplitz in 1924 \cite{toeplitz} gave an explicit construction of almost periodic functions (in the Bohr sense) on $\mathbb{R}$, where he as a device used arithmetic
progressions.) Toeplitz flows are closely related to the well understood odometer systems, and may in a certain sense be seen as the simplest Cantor minimal systems beyond the odometer systems.
Yet, surprisingly, Toeplitz flows exhibit a richness of properties that the odometer systems do not have. As examples of this we mention that any Choquet simplex can be realized as the set of invariant
probability measures for some Toeplitz flow. Furthermore, for every $0\leq t<\infty$ there exists a Toeplitz flow (in fact, an uncountable family of pairwise non-isomorphic Toeplitz flows) whose
topological entropy is equal to $t$.
\medbreak

An entirely new approach to study Toeplitz flows, and, for that matter, Cantor minimal systems in general -- thereby providing powerful new tools -- came from an unsuspected source : 
non-commutative $C^\ast$-algebras. In fact, by studying the so-called $C^\ast$-crossed product associated to a Cantor minimal system a complete isomorphism invariant turns out to be a 
special ordered abelian group, a so-called dimension group. This group is the $K^0$-group of the crossed product and it comes with a natural ordering. It turns out that the $K^0$-group can
be defined purely in dynamical terms, and that it is a complete invariant for (strong) orbit equivalence. This invariant is completely independent
of other invariants traditionally used to study dynamical systems, like spectral invariants and entropy. Furthermore, the invariant is in many cases effectively computable. We mention two papers
that illustrate this. The first is the paper by Durand, Host and Skau \cite{DHS} where substitutional dynamical systems are studied using the $K^0$-approach. The other is a paper by Gjerde
and Johansen \cite{gjerde-johansen} where Toeplitz flows are treated. In the latter paper a ``clean'' conceptual proof -- as seen from the theory of dimension groups -- is given of the Choquet simplex 
realization result alluded to above. (The original proof of that result is due to Downarowicz \cite{D}, after preliminary results had been obtained by Williams \cite{williams}.) Furthermore, 
the dimension groups (i.e.\ the $K^0$-groups), associated to Toeplitz flows are described in their paper in terms of Bratteli diagrams with what they call the EPN-property (which is the same as the ERS-property --
the term we prefer to use). In this paper we will both extend their results considerably, and we will give a more satisfactory description of the $K^0$-groups. In fact, the main thrust of this paper is to give the definite
characterization of the $K^0$-groups associated to Toeplitz flows. This is done in an intrinsic way, meaning that we give a characterization in terms of the group itself, while the Bratteli diagram realizations of the groups
in question play an important, but auxiliary role. By our characterization we can easily exhibit concrete examples of such groups, some of these will be described in Section \ref{sec:examples}. Furthermore, our
results underscore in a striking manner that the $K^0$-group of a Toeplitz flow $(X,T)$ does not ``see'' the entropy $h(T)$ of $T$, these two entities are independent. (See the Remark after Theorem \ref{th:real_simplex}.)
\medbreak

We remark that Toeplitz flows in contrast to substitution minimal systems are ``unstable''. By that we mean that whereas an induced system of a substitution system is again a substitution system, and a 
Cantor factor of a substitution system is either again a substitution system or an odometer, this it not true for Toeplitz flows in general. (Cf.\ \cite{DD}.)There is a ``tiny'' overlap between Toeplitz
flows and substitution minimal systems. In fact, some -- but not all -- substitution minimal systems associated to (primitive, aperiodic) substitutions of constant length are Toeplitz flows. (Cf.\ 
\cite[Theorem 6.03]{M}.)

\section{Main results}
\label{sec:results}
We formulate our main results, referring to Section \ref{sec:def} for definitions of some of the terms occurring in the statements of the theorems below. 

\begin{theorem}
\label{th:twoequalsets}
Let $0\leq t< \infty$. The following two sets are equal (up to order isomorphisms):
\begin{enumerate}[(i)]
 \item $\set{(K^0(X,T),K^0(X,T)^+}{(X,T) \text{ Toeplitz flow, } h(T)=t}$
 \item $\left\lbrace(G,G^+) \big| G \text{ simple dimension group containing a noncyclic subgroup }\right.$ \\
$ \left.\qquad \qquad\, \, H \text{ of rank one such that } H\cap G^+ \neq \{0\}\right\rbrace$
\end{enumerate}
\end{theorem}

\begin{theorem}
\label{th:realizing}
Let $(G,G^+,u)$ be a simple dimension group with order unit $u$, and assume the rational (and hence rank one) subgroup $\mathbb{Q}(G,u)$ is noncyclic. Let $0\leq t<\infty$. There exists a Toeplitz flow $(X,T)$ such that 
\begin{enumerate}[(i)]
 \item $(G,G^+,u)\cong (K^0(X,T),K^0(X,T)^+,\mathds{1})$
 \item The entropy $h(T)$ of $T$ equals $t$.
 \item The set of (continuous) eigenvalues of $T$ is $\left\lbrace\mathrm{e}^{2\pi\mathrm{i}s} \big| s\in \mathbb{Q}(G,u)\right\rbrace$
\end{enumerate}
\end{theorem}

\begin{theorem} 
\label{th:threeequalsets}
Let $0\leq t< \infty$. The following three sets of simple dimension groups with order units are equal (up to order isomorphisms preserving order units):
\begin{align*}
 \mathcal{T}_t &= \set{(K^0(X,T),K^0(X,T)^+,\mathds{1})}{(X,T) \text{ Toeplitz flow, } h(T)=t}\\
 \mathcal{G} &= \set{(G,G^+,u)}{G \text{ simple dimension group with } \mathbb{Q}(G,u) \text{ noncyclic}}\\
 \mathcal{B} &= \set{(K_0(V,E), K_0(V,E)^+,[\mathbf{1}])}{(V,E) \text{ simple Bratteli diagram with the ERS property}}.
\end{align*}
\end{theorem}

\begin{theorem}
\label{th:real_simplex}
Let $0\leq t< \infty$. For any Choquet simplex $K$ there exists an uncountable family of pairwise non-orbit equivalent Toeplitz flows
$\mathcal{A} = \{(X,T)\}$, such that for each $(X,T)\in\mathcal{A}$:
\begin{enumerate}[(i)]
 \item $K\cong M(X,T)$, i.e.\ $K$ is affinely homeomorphic to the $T$-invariant probability measures on $X$.
 \item $h(T)=t$.
\end{enumerate}

Now let $(Y,S)$ be any odometer and let $K$ be any Choquet simplex. There exists an uncountable family of pairwise non-strong orbit equivalent Toeplitz flows 
$\tilde{\mathcal{A}} = {(\tilde{X}, \tilde{T})}$, such that for each $(\tilde{X}, \tilde{T})\in \tilde{\mathcal{A}}$:
\begin{enumerate}[(i)]
 \item $K\cong M(\tilde{X}, \tilde{T})$
 \item $h(\tilde{T}) = t$
 \item The maximal equicontinuous factor of $(\tilde{X}, \tilde{T})$ is $(Y,S)$.
\end{enumerate}
\end{theorem}

\textit{Remark.}
As we alluded to above, Theorem \ref{th:real_simplex} yields in a striking way as a corollary that for Toeplitz flows orbit structure and entropy are independent entities. 
\medbreak

\section{Basic concepts and definitions and key background results.}
\label{sec:def}
\subsection{Dynamical systems}

Throughout this paper we will use the term dynamical system to mean a compact metric space $X$ together with a homeomorphism $\deff{T}{X}{X}$, and we will denote this by $(X,T)$. This induces in a natural way a $\mathbb{Z}$-action
on $X$. The orbit of $x\in X$ under this action is $\set{T^n x}{n\in \mathbb{Z}}$ and will be denoted by $\mathrm{orbit}_T(x)$. If all the orbits are dense in $X$ 
we say that $(X,T)$ is a \emph{minimal system}. It is  a simple observation that $(X,T)$ is minimal if and only if $TA=A$ for some closed
$A\subseteq X$ implies that $A=X$ or $\varnothing$.

(We will denote the natural numbers $\{1,2,3,\dots\}$ by $\mathbb{N}$, the integers by $\mathbb{Z}$, the rational numbers by $\mathbb{Q}$, the real numbers by $\mathbb{R}$. Also, let 
$\mathbb{Z}^+ = \{0,1,2,\dots\}$, $\mathbb{Q}^+ = \{r\in\mathbb{Q}\,|\, r\geq 0\}$, $\mathbb{R}^+ = \{t\in\mathbb{R}\,|\, t\geq 0\}$.)

\begin{definition}
We say that a dynamical system $(Y,S)$ is a \emph{factor} of $(X,T)$ and that $(X,T)$ is an \emph{extension} of $(Y,S)$ if there exists a continuous surjection $\deff{\pi}{X}{Y}$ which 
satisfies $S(\pi(x))=\pi(Tx),\,\forall x\in X$. We call $\pi$ a \emph{factor map}.
If $\pi$ is a bijection then we say that $(X,T)$ and $(Y,S)$ are \emph{conjugate}, and we write $(X,T)\cong (Y,S)$. We say that $(X,T)$ is \emph{flip conjugate} to $(Y,S)$ if $(X,T)\cong (Y,S)$ or $(X,T)\cong (Y,S^{-1})$.
\end{definition}

\begin{definition}
\label{def:orbit}
The dynamical systems $(X,T)$ and $(Y,S)$ are \emph{orbit equivalent} if there exists a homeomorphism $\deff{F}{X}{Y}$ such that $F(\mathrm{orbit}_T(x)) = \mathrm{orbit}_S(F(x))$ for all $x\in X$. We call $F$ an 
\emph{orbit map}.
\end{definition}

\proc{Remark.}
Clearly flip conjugacy implies orbit equivalence. It is a fact that if $X$ (and hence $Y$) is a connected space then orbit equivalence between $(X,T)$ and $(Y,S)$ implies flip conjugacy. This has the consequence that the study of
orbit equivalence is only interesting as it pertains to \emph{Cantor minimal systems} $(X,T)$, i.e.\ $X$ is a Cantor set on which $T$ acts minimally. The $K$-theoretic invariant we are going to introduce is an invariant
for orbit equivalence, and so we will assume henceforth that our dynamical systems are Cantor minimal, even though some of the subsequent definitions apply to more general systems. 
\medbreak

Let $(X,T)$, $(Y,S)$ and $F$ be as in Definition \ref{def:orbit}, both $(X,T)$ and $(Y,S)$ being Cantor minimal systems. For each $x\in X$ there exists a unique integer $n(x)$ (respectively, $m(x)$)  such that $F(Tx)= S^{n(x)}(F(x))$,
$F(T^{m(x)}x) = S(F(x))$. We call $\deff{m,n}{X}{\mathbb{Z}}$ the orbit cocycles associated to the orbit map $F$.

\begin{definition}
 We say that $(X,T)$ and $(Y,S)$ are strong orbit equivalent if there exists an orbit map $\deff{F}{X}{Y}$ such that each of the two associated orbit cocycles $\deff{m,n}{X}{\mathbb{Z}}$ have at 
most one point of discontinuity.
\end{definition}

\proc{Remark.}
Boyle proved that if the orbit cocycles are continuous for all $x\in X$, then $(X,T)$ and $(Y,S)$ are flip conjugate. So strong orbit equivalence is in a sense the mildest weakening possible of flip conjugacy. 
\medbreak

We shall need the concept of \emph{induced transformation} and \emph{Kakutani equivalence}.

\begin{definition}
\label{def:induced}
 Let $(X,T)$ be a Cantor minimal system and let $A$ be a clopen subset of $X$ (hence $A$ is again a Cantor set). Let $\deff{T_A}{A}{A}$ be the \emph{first return map}, i.e.\ $T_A(x)=T^{r_A(x)} x$, where 
$r_A(x) = \min\set{n\in \mathbb{N}}{T_A^nx\in A}$. We say that $(A,T_A)$, which is again Cantor minimal, is the induced system of $(X,T)$ with respect to $A$. 
\end{definition}

\begin{definition}
\label{def:kaku}
 The Cantor minimal systems $(X,T)$ and $(Y,S)$ are \emph{Kakutani equivalent} if (up to conjugacy) they have a common induced system. (One can show that this is an equivalence relation on the family of Cantor minimal systems.)
\end{definition}

\begin{definition}
 Let $(X,T)$ and $(Y,S)$ be (Cantor) minimal systems, $(Y,S)$ being a factor of $(X,T)$ by a map $\deff{\pi}{X}{Y}$. If there exists a point $x\in X$ such that $\pi^{-1}\left(\pi(x)\right) = \{x\}$ we say that $\pi$ is
almost one-to-one and we say that $(X,T)$ is an almost one-to-one extension of $(Y,S)$. (One can show that the set consisting of points $x$ in $X$ satisfying the above condition is a dense $G_\delta$-set of $X$.)  
\end{definition}

\subsection{Toeplitz flows}
\label{sec:toeplitz}
\begin{definition}
 $(X,T)$ is \emph{expansive} if there exists $\delta>0$ such that if $x\neq y$ then $\mathrm{sup}_n \mathrm{d}(T^nx, T^n y)>\delta$, where $\mathrm{d}$ 
is a metric that gives the topology of $X$. (Expansiveness is independent of the metric as long as the metric gives the topology of $X$.)
\end{definition}

Let $\Lambda=\{a_1,a_2,\dots,a_n\},\, n\geq 2$, be a finite alphabet and let $Z=\Lambda^\mathbb{Z}$ be the set of all bi-infinite sequences of symbols from
$\Lambda$ with $Z$ given the product topology -- thus $Z$ is a Cantor set. Let $\deff{S}{Z}{Z}$ denote the shift map, $\deff{S}{(x_n)}{(x_{n+1})}$.
If $X$ is a closed subset of $Z$ such that $S(X)=X$, we say that $(X,S)$ is a \emph{subshift}, where we denote the restriction of $S$ to $X$ again by $S$.
Subshifts are easily seen to be expansive. We state the following well-known fact as a proposition. (Cf.\ \cite[Theorem 5.24]{walters}.)

\begin{proposition}
 \label{prop:subshift}
Let $(X,T)$ be a Cantor minimal system. Then $(X,T)$ is conjugate to a minimal subshift on a finite alphabet if and only if $(X,T)$ is expansive.
\end{proposition}

As a general reference on Toeplitz flows we refer to \cite{williams}. (Cf.\ also \cite{D}.)

\begin{definition}
 Let $\eta=(\eta(n))_{n\in\mathbb{Z}}\in\Lambda^\mathbb{Z}$, where $\Lambda$ is a finite alphabet. Then we define for $\sigma\in \Lambda$, $p\in\mathbb{N}$
\[ \per{p}{\eta,\sigma} = \set{n\in\mathbb{Z}}{\eta(n+mp)=\sigma,\, \forall m\in\mathbb{Z}}.\]
Let
\[\per{p}{\eta}= \bigcup_{\sigma \in \Lambda} \per{p}{\eta,\sigma}.\]
We say that $\eta$ is a \emph{Toeplitz sequence} if $\,\bigcup_{p\in\mathbb{N}}\per{p}{\eta} = \mathbb{Z}$.
\end{definition}

By the \emph{$p$-skeleton} of $\eta$ we will mean the part of $\eta$ which is periodic with period $p$; more precisely, we define the $p$-skeleton to be the sequence obtained from $\eta$ by replacing $\eta(n)$ with a new symbol
$\ast$ for all $n\not\in \per{p}{\eta}$. We say that $p$ is an \emph{essential period} of $\eta$ if the $p$-skeleton of $\eta$ is not periodic
with any smaller period. The least common multiple of two essential periods is again an essential period, a fact which is easily verified. 

\begin{definition}
 Assume that a Toeplitz sequence $\eta$ is non-periodic. A \emph{periodic structure} for $\eta$ is a strictly increasing sequence $(p_i)_{i\in\mathbb{N}}$
such that $p_i$ is an essential period of $\eta$ for all $i$, $p_i | p_{i+1}$ and $\bigcup_{i=1}^\infty \per{p_i}{\eta}=\mathbb{Z}$. (It is a fact that
a periodic structure always exists for a (non-periodic) Toeplitz sequence.)
\end{definition}

\begin{definition}
 Let $\eta\in\Lambda^\mathbb{Z}$ be a Toeplitz sequence. The dynamical system $(\bar{\mathcal{O}(\eta)}, S)$ is called a \emph{Toeplitz flow}, 
where $\mathcal{O}(\eta) = \mathrm{orbit}_S(\eta)$ and $S$ is the shift map. 
\end{definition}

Periodic sequences are Toeplitz sequences. Every Toeplitz sequence is almost periodic, i.e.\ every word occurring in $\eta$ appears with bounded gap 
between successive occurrences. Hence $(\bar{\mathcal{O}(\eta)}, S)$ is minimal. In the sequel we will only consider non-periodic Toeplitz flows, and 
so Toeplitz flows are expansive Cantor minimal systems. 

Let $(p_i)_{i\in\mathbb{N}}$ be a periodic structure for the (non-periodic) Toeplitz sequence $\eta$, and denote the associated Toeplitz flow by 
$(Y,S)$. Let $(G_\mathfrak{a},\rho_{\hat{\mathds{1}}})$ denote the \emph{odometer} (also called \emph{adding machine}) associated to the $\mathfrak{a}$-adic group
\[G_\mathfrak{a} = \prod_{i=1}^\infty \left\{0,1,\dots \frac{p_i}{p_{i-1}}-1\right\},\]
where $\mathfrak{a}= \left\{\frac{p_i}{p_{i-1}}\right\}_{i\in\mathbb{N}}$ (we set $p_0=1$) and where $\rho_{\hat{\mathds{1}}}(x) = x+\hat{1}$, where $\hat{1} = (1,0,0,\dots)$.
We note that $G_\mathfrak{a}$ is naturally isomorphic to the inverse limit group  
\[\quotient{\mathbb{Z}}{p_1\mathbb{Z}} \stackrel{\phi_1}{\longleftarrow}\quotient{\mathbb{Z}}{p_2\mathbb{Z}} \stackrel{\phi_2}{\longleftarrow} \quotient{\mathbb{Z}}{p_3\mathbb{Z}} \stackrel{\phi_3}{\longleftarrow} \cdots\]
where $\phi_i(n)$ is the residue of $n$ modulo $p_i$. It is a fact that the family consisting of compact groups $G$ that are both monothetic (i.e.\ contains a dense copy
of $\mathbb{Z}$ , which of course implies that $G$ is abelian) and Cantor (as a topological space), coincides with the family of 
$\mathfrak{a}$-adic groups. It is also noteworthy that all minimal rotations (in particular rotations by $\hat{1}$) on such groups are conjugate. This is a consequence of the fact that the 
dual group of an $\mathfrak{a}$-adic group is a torsion group.
If $\mathfrak{a} = \{p\}_{i\in\mathbb{N}}$, where $p$ is a prime, then $G_\mathfrak{a}$ is the $p$-adic integers. (We refer to \cite[Vol 1]{HR} for background information on $\mathfrak{a}$-adic groups.)
\medbreak

Recall that every (minimal) dynamical system has a \emph{maximal equicontinuous factor}, the latter being, by a well-known theorem, conjugate to a minimal rotation on a compact abelian group $G$, and so, in particular, is uniquely ergodic.
Specifically, let $(G,\rho_g)$ be the maximal equicontinuous factor of $(X,T)$, where $\deff{\rho_g}{G}{G}$ is rotation by $g$, i.e.\ $\rho_g(x) = x+g$ for $x\in G$, and $\deff{\Phi}{X}{G}$ is the factor map. Let $(H,\rho_h)$ be another minimal
group rotation factor of $(X,T)$  where $\deff{\Psi}{X}{H}$ is the factor map. Then $(H,\rho_h)$ is a factor of $(G,\rho_g)$ by a factor map $\deff{\pi}{G}{H}$ such that $\Psi=\pi\circ\Phi$.

One can detect the maximal equicontinuous factor $(G,\rho_g)$ of $(X,T)$ by determining the (continuous) eigenvalues $\Gamma$ of $T$, i.e.\ 
\[\Gamma = \set{\lambda \in\mathbb{T}}{f\circ T = \lambda f \text{ for some } 0\neq f\in \mathcal{C}(X)},\]
where $\mathbb{T}$ is the unit circle. Now $\Gamma$ is a countable and discrete subgroup of the discrete circle group $\mathbb{T}_d$, and $G= \hat{\Gamma}$ is the compact dual group with $g\in G$ being the character on $\Gamma$ defined
by $g(\gamma) = \gamma$, $\forall \gamma \in \Gamma$.  Incidentally, the dual group $\hat{G_\mathfrak{a}}$ of the odometer group $G_\mathfrak{a}$ has only torsion elements and is equal to the following subgroup of $\mathbb{T}_d$,
namely 
\[\hat{G_\mathfrak{a}} =  \left\lbrace \mathrm{e}^{\tfrac{2\pi\mathrm{i}l}{a_1a_2\cdots a_k}}\,\Big|\,\mathfrak{a} = \{a_i\}_{i\in\mathbb{N}}, k\in\mathbb{N}, l\in\mathbb{Z}\right\rbrace.\]

If $(X,T)$ is a Toeplitz flow with a periodic structure $(p_i)_{i\in\mathbb{N}}$, then the maximal equicontinuous factor is the odometer system associated to $\mathfrak{a}=\left(\frac{p_i}{p_{i-1}}\right)_{i\in\mathbb{N}}$.
Furthermore, the factor map $\deff{\pi}{X}{G_\mathfrak{a}}$ is an almost one-to-one map. The converse is also true. We state two theorems that establish this fact. 

\begin{theorem}[\cite{P}]
 \label{th:paul}
Let $(X,T)$ be a minimal almost one-to-one extension of a (minimal) group rotation system $(G,\rho_g)$. Then $(G,\rho_g)$ is the maximal equicontinuous factor of $(X,T)$.
\end{theorem}

\begin{theorem}[\cite{MP}]
 \label{th:markley-paul}
An expansive Cantor minimal system $(X,T)$ is a Toeplitz flow if and only if $(X,T)$ is an almost one-to-one extension of an odometer system $(G_\mathfrak{a},\rho_{\hat{1}})$. Furthermore, $(G_\mathfrak{a},\rho_{\hat{1}})$ is 
the maximal equicontinuous factor of $(X,T)$.
\end{theorem}

A more detailed analysis of the situation is found in \cite{williams}. We summarize by listing some salient points from \cite{williams} which is relevant for this paper:
\medbreak
Let $\eta$ be a (non-periodic) Toeplitz sequence and let $(\bar{\mathcal{O}(\eta)}, S)$ be the associated Toeplitz flow. Let $(p_i)_{i\in\mathbb{N}}$ be a periodic structure for $\eta$. Define $A_n^i = \set{S^m\eta}{m\equiv n (\mod p_i)}$.
\medbreak
\begin{enumerate}
 \item \label{pt:skeleton}$\bar{A_n^i}$ is the set of all $\omega\in\bar{\mathcal{O}(\eta)}$ with the same $p_i$-skeleton as $S^n(\eta)$.
\medbreak
\item \label{pt:partition}$\set{\bar{A_n^i}}{n\in\quotient{\mathbb{Z}}{p_i\mathbb{Z}}}$ is a partition of $\bar{\mathcal{O}(\eta)}$ into clopen sets.
\medbreak
\item \label{pt:nesting} $\bar{A_n^i} \supset \bar{A_m^j}$ for $i<j$ and $m\equiv n (\mod{p_i})$.
\medbreak
\item \label{pt:shift} $S\bar{A_n^i}=\bar{A_{n+1}^i}$.
\medbreak
\item \label{pt:} $\omega, \theta \in \bigcap_{i=1}^\infty \bar{A_{n_i}}$ (where $n_i\equiv n_j (\mod p_i)$ for $j\geq i$) if and only if $\omega$ and $\theta$ have the same $p_i$-skeleton for all $i\in\mathbb{N}$. In particular, 
$\bigcap_{i=1}^\infty \bar{A_{n_i}} = \{\omega\}$ if and only if $\omega$ is a Toeplitz sequence. This implies that $\pi^{-1}(\pi(\omega)) = \{\omega\}$ if and only if $\omega$ is a Toeplitz sequence. Here $\deff{\pi}{\bar{\mathcal{O}(\eta)}}{G_\mathfrak{a}}$
is the factor map, $(G_\mathfrak{a},\rho_{\hat{1}})$ being the maximal equicontinuous factor of $(\bar{\mathcal{O}(\eta)}, S)$.
\end{enumerate}
\medbreak
Since $\mathrm{Per}_{p_i}(\eta)$ is periodic it has density in $\mathbb{Z}$ given by
\[d_i = \frac{1}{p_i} \left|\set{n\in\quotient{\mathbb{Z}}{p_i\mathbb{Z}}}{n\in\mathrm{Per}_{p_i}(\eta)}\right|\]
where $|E|$ denotes the number of elements in the set $E$. The $d_i$ are increasing, and we set $d=\lim_{i\rightarrow \infty} d_i$.

\begin{definition}
\label{def:regular}
 The Toeplitz sequence $\eta$ is \emph{regular} if $d=1$.
\end{definition}

\proc{Notation.}
Let $(X,T)$ be a dynamical system. We denote by $M(X,T)$ the set of $T$-invariant probability measures on $X$. We say that $(X,T)$ is \emph{uniquely ergodic} if $M(X,T)$ is a singleton.
\medbreak

\proc{Fact.}
If $(X,T)$ is a minimal dynamical system, then $M(X,T)$ is a Choquet simplex (where $M(X,T)$ is given the $w^\ast$-topology).
\medbreak

\begin{theorem}[\cite{JK}]
 If $(\bar{\mathcal{O}(\eta)}, S)$ is a regular Toeplitz flow (i.e.\ $\eta$ is a regular Toeplitz sequence), then it is uniquely ergodic and the (topological) entropy $h(S)$ is zero.
\end{theorem}

We will return to Toeplitz flows and their properties after we have introduced the (ordered) K-theory associated to such systems.

\subsection{Bratteli diagrams and dimension groups.}
(As general references for the material in this section we refer to \cite{effros}, \cite{HPS} and \cite{GPS}.)
\subsubsection{Bratteli diagrams} 
\begin{figure}
\centering
 \includegraphics{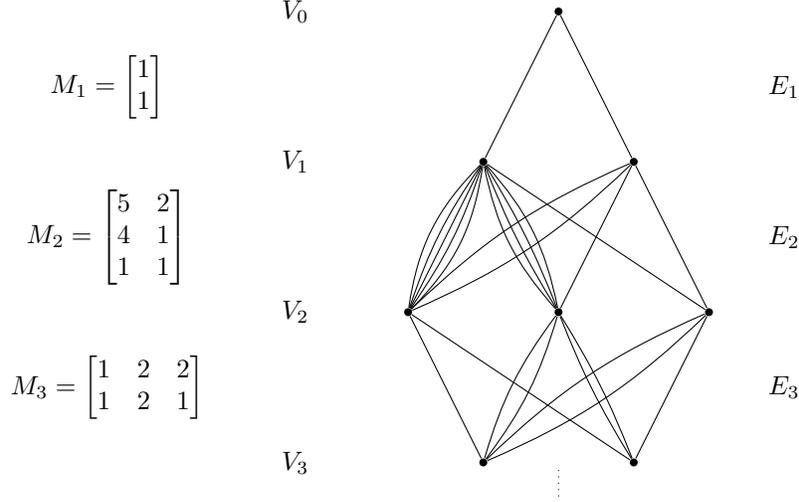}
\caption{An example of a Bratteli diagram}
\label{fig:Bd}
\end{figure}

A Bratteli diagram $(V,E)$ consists of a set of vertices $V= \sqcup_{n=0}^\infty V_n$ and a set of edges $E=\sqcup_{n=1}^\infty E_n$, where the $V_n$'s and the $E_n$'s are finite disjoint sets and where
$V_0=\{v_0\}$ is a one-point set. The edges in $E_n$ connect vertices in $V_{n-1}$ with vertices in $V_n$. If $e$ connects $v \in V_{n-1}$ with $u \in V_n$ we write $s(e)=v$ and $r(e)=u$, where $\deff{s}{E_n}{V_{n-1}}$ and
$\deff{r}{E_n}{V_n}$ are the source and range maps, respectively. We will assume that $s^{-1}(v)\neq \varnothing$ for all $v\in V$ and that $r^{-1}(v)\neq \varnothing$ for all $v\in V\backslash V_0$. A Bratteli diagram
can be given a diagrammatic presentation with $V_n$ the vertices at level $n$ and $E_n$ the edges between $V_{n-1}$ and $V_n$. If $|V_{n-1}|=t_{n-1}$ and $|V_{n}|=t_{n}$ then the edge set $E_n$ is described by
a $t_n\times t_{n-1}$ incidence matrix $M_n=(m_{ij}^n)$, where $m_{ij}^n$ is the number of edges connecting $v_i^n\in V_n$ with $v_j^{n-1}\in V_{n-1}$ (see Figure \ref{fig:Bd}). Let $k,l\in \mathbb{Z}^+$ with $k<l$ and let $E_{k+1}\circ E_k\circ \cdots
\circ E_l$ denote all the paths from $V_k$ to $V_l$.  Specifically, 
\[E_{k+1}\circ E_k\circ \cdots \circ E_l = \set{(e_{k+1},\cdots,e_l)}{e_i\in E_i, i=k+1,\dots,l ; r(e_i)=s(e_{i+1}),i=k+1,\dots,l-1}.\]
We define $r\left((e_{k+1},\cdots,e_l)\right) = r(e_l)$ and $s\left((e_{k+1},\cdots,e_l)\right) = s(e_{k+1})$. Notice that the corresponding incidence matrix is the product $M_{l}M_{l-1}\cdots M_{k+1}$ of the incidence matrices.

\begin{definition}
\label{def:ers}
 The Bratteli diagram $(V,E)$ with incidence matrices $(M_n)_{n=1}^\infty$ has the ERS-property (ERS = Equal Row Sum) if the row sums of the incidence matrices are constant. Let the constant row sum of $M_n$ be $r_n$. We associate
the supernatural number $\prod_{n=1}^\infty r_n$ to $(V,E)$. (See the comments after Definition \ref{def:rational}.)
\end{definition}

\begin{definition}
Given a Bratteli diagram $(V,E)$ and a sequence $0 = m_0 < m_1 < m_2 < \cdots$ in $\mathbb{Z^+}$, we define the \emph{telescoping} of $(V,E)$ to $\{m_n\}$ as $(V',E')$, where $V_n'=V_{m_n}$ and 
$E_n' = E_{m_{n-1}+1}\circ \cdots \circ E_{m_n}$, and the source and the range maps are as above.
\end{definition}

\proc{Remark.}
Observe that the ERS-property is preserved under telescoping, and that is also the case for the associated supernatural number.

\begin{definition}
We say that the Bratteli diagram $(V,E)$ is \emph{simple} if there exists a telescoping of $(V,E)$ such that the resulting Bratteli diagram $(V',E')$ has full connection between all consecutive levels, i.e.\ the entries of all the 
incidence matrices are non-zero.
\end{definition}

Given a Bratteli diagram $(V,E)$ we define the infinite path space associated to $(V,E)$, namely 
\[X_{(V,E)} = \set{(e_1, e_2, \dots)}{e_i\in E_i, r(e_i)=s(e_{i+1}); \quad \forall i \geq 1}.\]
Clearly $X_{(V,E)} \subseteq \prod_{n=1}^\infty E_n$, and we give $X_{(V,E)}$ the relative topology, $\prod_{n=1}^\infty E_n$ having the product topology. Loosely speaking this means that two paths in $X_{(V,E)}$ are close
if the initial parts of the two paths agree on a long initial stretch. Also, $X_{(V,E)}$ is a closed subset of $\prod_{n=1}^\infty E_n$, and is compact. 

\medbreak
Let $p=(e_1, e_2, \dots, e_n)\in E_1\circ\cdots\circ E_n$ be a finite path starting at $v_0\in V_0$. We define the \emph{cylinder set} $U(p)=\set{(f_1,f_2,\dots)\in X_{(V,E)}}{f_i = e_i, i = 1,2,\dots,n}$. The collection of 
cylinder sets is a basis for the topology on $X_{(V,E)}$. The cylinder sets are clopen sets, and so $X_{(V,E)}$ is a compact, totally disconnected metric space -- metric because the collection of cylinder sets is countable. 
If $(V,E)$ is simple then $X_{(V,E)}$ has no isolated points, and so $X_{(V,E)}$ is a Cantor set. (We will in the sequel disregard the trivial case where $|X_{(V,E)}|$ is finite.)

Let $P_n =E_1\circ\cdots\circ E_n$ be the set of finite paths of length $n$ (starting at the top vertex). We define the truncation map $\deff{\tau_n}{X_{(V,E)}}{P_n}$ by $\tau_n\left((e_1,e_2,\dots)\right) = (e_1,e_2,\dots,e_n)$.
If $m\geq n$ we have the obvious truncation map $\deff{\tau_{m,n}}{P_m}{P_n}$.

There is an obvious notion of isomorphism between Bratteli diagrams $(V,E)$ and $(V',E')$; namely, there exists a pair of bijections between $V$ and $V'$ preserving the gradings and intertwining the respective source and range
maps. Let $\sim$ denote the equivalence relation on Bratteli diagrams generated by isomorphism and telescoping. One can show that $(V,E)\sim(V',E')$ iff there exists a Bratteli diagram $(W,F)$ such that telescoping $(W,F)$ to odd
levels $0<1<3<\cdots$ yields a diagram isomorphic to some telescoping of $(V,E)$, and telescoping $(W,F)$ to even levels $0<2<4<\cdots$ yields a diagram isomorphic to some telescoping of $(V',E')$.

\subsubsection{Dimension groups}
By an ordered group we shall mean a countable abelian group $G$ together with a subset $G^+$, called the positive cone, such that 
\begin{enumerate}
 \item $G^+-G^+ = G$
 \item $G^+\cap (-G^+) = \{0\}$
 \item $G^+ + G^+ \subset G^+$
\end{enumerate}
We shall write $a\leq b$ if $b-a\in G^+$. We say that an ordered group is \emph{unperforated} if $a\in G$ and $na\in G^+$ for some $a\in G$ and $n\in\mathbb{N}$ implies that $a\in G^+$. Observe that an unperforated group
is torsion free. By an \emph{order unit} for $(G,G^+)$ we mean an element $u\in G^+$ such that for every $a\in G$, $a\leq nu$ for some $n\in \mathbb{N}$. 

\begin{definition}
A \emph{dimension group} $(G,G^+,u)$ with distinguished order unit $u$ is an unperforated ordered group $(G,G^+)$ satisfying the \emph{Riesz interpolation property}, i.e.\ given $a_1, a_2, b_1, b_2 \in G$ with $a_i\leq b_j$ ($i,j=1,2$),
there exists $c\in G$ with $a_i\leq  c \leq b_j$ $(i,j=1,2)$.
\end{definition}

We write $(G_1,G_1^+,u_1)\cong (G_2,G_2^+,u_2)$ if there exists an order isomorphism $\deff{\phi}{G_1}{G_2}$, i.e.\ $\phi$ is a group isomorphism such that $\phi(G_1^+)=G_2^+$, and
$\phi(u_1)=u_2$.

\medbreak

To a Bratteli diagram $(V,E)$ we can associate an ordered group, which we will denote by $K_0(V,E)$ (because of its connection to (ordered) K-theory). Let $V=\sqcup_{n=0}^\infty V_n$ and let $(M_n)_{n=1}^\infty$ be
the incidence matrices. Then we have a system of simplicially ordered groups and positive maps 
\[(\mathbb{Z} = ) \,\mathbb{Z}^{\abs{V_0}}\stackrel{M_1}{\longrightarrow}\mathbb{Z}^{\abs{V_1}}\stackrel{M_2}{\longrightarrow}\mathbb{Z}^{\abs{V_2}}\longrightarrow \cdots\]
where the positive homomorphism $\deff{M_n}{\mathbb{Z}^{|V_{n-1}|}}{\mathbb{Z}^{|V_{n}|}}$ is given by matrix multiplication with the incidence matrix $M_n$. ($\mathbb{Z}^{|V_{n}|}$ is a column 
vector, and an element in $\mathbb{Z}^{|V_{n}|}$ is positive if all its entries are non-negative.) By definition $K_0(V,E)$ is the inductive limit of the system above, and $K_0(V,E)$ is given the induced order. We denote the positive
cone by $K_0(V,E)^+$. $K_0(V,E)$ has a distinguished \emph{order unit}, namely the element $[\mathbf{1}]$ in $K_0(V,E)^+$ corresponding to the element $1\in \mathbb{Z}^{|V_0|}= \mathbb{Z}$. The triple $(K_0(V,E), K_0(V,E)^+, [\mathbf{1}])$
denotes the countable, ordered abelian group $K_0(V,E)$ with positive cone $K_0(V,E)^+$ and distinguished order unit $[\mathbf{1}]$. We will sometimes in the sequel for short only write $K_0(V,E)$, the ordering and the order unit being implicitly
understood. 
\medbreak

One can show that $(V,E)\cong(V',E')$ if and only if $K_0(V,E)$ is order isomorphic to $K_0(V',E')$ by a map sending the distinguished order unit of $K_0(V,E)$ to the distinguished order unit of $K_0(V',E')$.

It is straightforward to verify that $(K_0(V,E), K_0(V,E)^+)$ is a dimension group. The converse however, is not obvious and we state it as a theorem formulated in such a way that it suits our purpose. 

\begin{theorem}[\cite{EHS}]
 Let $(G,G^+,u)$ be a dimension group with distinguished order unit $u$. Then there exists a Bratteli diagram $(V,E)$ such that $(G,G^+,u)\cong (K_0(V,E),K_0(V,E)^+,[\mathbf{1}])$.
\end{theorem}

\medbreak
A dimension group $(G,G^+)$ is \emph{simple} if it contains no non-trivial order ideals. An order ideal is a subgroup $J$ such that $J= J^+-J^+$ (where $J^+ = J\cap G^+$) and $0\leq a\leq b \in J$ implies $a\in J$.
It is easily seen that $(G,G^+)$ is a simple dimension group if and only if every $a\in G^+\backslash \{0\}$ is an order unit. In the sequel we will exclusively work with noncyclic (i.e.\ $G\not\cong \mathbb{Z}$) simple dimension groups $(G,G^+)$.
\medbreak

Let $(G,G^+,u)$ be a simple dimension group with distinguished order unit $u$. We say that a homomorphism $\deff{p}{G}{\mathbb{R}}$ is a \emph{state} if $p$ is positive (i.e.\ $p(G^+)\geq 0$) and $p(u)=1$. Denote the
collection of all states on $(G,G^+,u)$ by $S_u(G)$. Now $S_u(G)$ is a convex compact subset of the locally convex space $\mathbb{R}^G$ with the product topology. In fact, one can show that $S_u(G)$ is a Choquet simplex. It
is a fact that $S_u(G)$ determines the order on $G$. In fact,
\[G^+ = \set{a\in G}{p(a)>0, \forall p\in S_u(G)}\cup \{0\}.\]

\begin{definition}
Let $(G,G^+)$ be a simple dimension group and let $u\in G^+\backslash\{0\}$. We say that $a\in G$ is \emph{infinitesimal} if $-\epsilon u \leq a\leq \epsilon u$ for all $0<\epsilon\in \mathbb{Q}^+$. (If
$\epsilon = \frac{p}{q},\, p,q\in \mathbb{N}$, then  $a\leq \epsilon u$ means that $qa \leq pu$. It is easy to see that the definition does not depend upon the particular order unit $u$.) An equivalent definition is: $a\in G$
is infinitesimal if $p(a)=0$ for all $p\in S_u(G)$. The collection of infinitesimal elements of $G$ form a  subgroup, \emph{the infinitesimal subgroup of $G$}, which we denote by $\mathrm{Inf}(G)$.
\end{definition}

\proc{Remark}
The quotient group $\quotient{G}{\mathrm{Inf}(G)}$ is again a simple dimension group in the induced order, and the infinitesimal subgroup of $\quotient{G}{\mathrm{Inf}(G)}$ is trivial. Also, an order unit for $G$ maps to an order unit for 
$\quotient{G}{\mathrm{Inf}(G)}$.
\medbreak

The following theorem summarizes some facts that are relevant for our situation, and the proof can be found in \cite{effros}.

\begin{theorem}
\label{th:simplex}
 Let $(G,G^+,u)$ be a simple (noncyclic) dimension group with distinguished order unit $u$.The map $\deff{\theta}{G}{\mathrm{Aff}(S_u(G))}$ from $G$ to the additive group of continuous affine functions on the Choquet simplex $S_u(G)$ defined
by $\theta(g)(p)=p(g)$ is a strict order preserving map (i.e.\ $\theta(g)(p)>0$ for all $p\in S_u(G)$ implies $g>0$). Furthermore, $\mathrm{Im}(\theta)$ is dense in $\mathrm{Aff}(S_u(G))$ in the uniform norm and contains
 the constant function $1$, and $\mathrm{ker}(\theta)=\mathrm{Inf}(G)$.

Conversely, suppose $K$ is a Choquet simplex and $H$ is a countable dense subgroup of $\mathrm{Aff}(K)$, and that $\deff{\theta}{G}{H}$ is a homomorphism of a torsion free countable abelian group $G$ onto $H$. Then letting 
\[G^+ = \set{g\in G}{\theta(g)(p)>0 \text{ all } p\in K}\cup \{0\}\]
we get that $(G,G^+)$ is a simple dimension group such that $\mathrm{Inf}G = \mathrm{ker}(\theta)$. In particular, if $G=H$ (with $\theta$ the identity map) and $G$ contains the constant function $1$, then $\mathrm{Inf}(G)=\{0\}$ and
$S_1(G)$ is affinely homeomorphic to $K$ by the map sending $k\in K$ to $\deff{\hat{k}}{G}{\mathbb{R}}$, where $\hat{k}(g) = g(k)$, $g\in G$.
\end{theorem}

\begin{definition}
\label{def:rational}
By a \emph{rational group} $H$ we shall mean a (additive) subgroup of $\mathbb{Q}$ that contains $\mathbb{Z}$. We say that $H$ is a \emph{noncyclic} rational group if $H$ is not isomorphic to $\mathbb{Z}$. Clearly $(H,H^+,1)$
is a simple dimension group with distinguished order unit $1$, where $H^+ = H\cap \mathbb{Q}^+$.
\end{definition}

Since rational dimension groups are going to play such an important role in this paper we will give a short description of them (Cf.\ \cite[Chapter XIII, Section 85]{fuchs}). First of all they are exactly the countable torsion-free groups of rank one. Let 
$n=p_1^{k_1}p_2^{k_2}p_3^{k_3}\cdots$ be a \emph{supernatural number}, where $p_1,p_2,p_3,\dots$ are the primes $2,3,5,\dots$ listed in increasing order, and $0\leq k_i \leq \infty$ for each $i$. (If $(r_n)_{n=1}^\infty$ is a sequence
of natural numbers, then we get a supernatural number $\prod_{n=1}^\infty r_n$ in an obvious way by factoring each $r_n$ into a product of primes.) Clearly $n\in\mathbb{N}$ if and 
only if $k_i< \infty$ for all $i$, and $k_i=0$ for all but finitely many $i$'s. If $m = p_1^{l_1}p_2^{l_2}p_3^{l_3}\cdots$ is another supernatural number we multiply $n$ and $m$ as $nm=p_1^{k_1+l_1}p_2^{k_2+l_2}p_3^{k_3+l_3}\cdots$,
where $k_i+l_i = \infty$ if either $k_i$ or $l_i$ are equal to $\infty$. We say that $m$ divides $n$ (notation $m|n$) if $l_i\leq k_i$ for all $i$. For $n=p_1^{k_1}p_2^{k_2}p_3^{k_3}\cdots$ let 
\[G(n) = \set{\frac{a}{b}}{a\in\mathbb{Z}, b\in\mathbb{N}, b|n}.\] 
Then $G(n)$ is a rational group and all rational groups are of this form. Furthermore, $G(n)$ is isomorphic to $G(m)$ if and only if there exists $a,b\in\mathbb{N}$ such that $an=bm$. In particular, all groups $G(n)$, where $n\in\mathbb{N}$, 
are isomorphic to $G(1) = \mathbb{Z}$. We note that $G(n)$ is $p$-divisible (i.e.\ for every $a\in G$ there exists $x\in E$ such that $px=a$) for some prime $p$ if and only if $p$ occurs with infinite multiplicity in the factorization of $n$. The group $G(n)$ can be made into a dimension group in exactly two ways, namely by letting the positive cone $G(n)^+$ be $G(n)\cap \mathbb{Q}^+$ or $-G(n)\cap \mathbb{Q}^+$, respectively.
If $\mathfrak{a} = (a_1,a_2,\dots)$ is a sequence of natural numbers with $a_i\geq 2$ for all $i$, we associate the noncyclic rational group $\left\lbrace \left.\frac{m}{a_1a_2\cdots a_k}\right| m\in \mathbb{Z}, k\in \mathbb{N}\right\rbrace$.
Clearly this is the same as the group $G(n)$, where $n$ is the  supernatural number $\prod_{i=1}^\infty a_i$. (Here the factorization of $\prod_{i=1}^\infty a_i$ into products of primes is obviously understood.)
\begin{definition}
\label{def:simple}
Let $(G,G^+,u)$ be a simple dimension group with order unit $u$. We define the rational subgroup of $G$, denoted $\mathbb{Q}(G,G^+,u)$ (or $\mathbb{Q}(G,u)$ for short), to be
\[\mathbb{Q}(G,u) = \set{g\in G}{ng=mu \text{ for some } n\in \mathbb{N}, m\in \mathbb{Z}}.\]
\end{definition}

\begin{proposition}
\label{prop:quotient}
 Let $\mathbb{Q}(G,u)$ be as in Definition \ref{def:simple} The map $\deff{\Gamma}{\mathbb{Q}(G,u)}{\mathbb{Q}}$ defined by $\Gamma(g)=\frac{m}{n}$ if $ng=mu$, is an injective order homomorphism sending $u$ to $1\in \mathbb{Z}$,
where $\mathbb{Q}(G,u)^+ = \mathbb{Q}(G,u)\cap G^+$. Thus $\mathbb{Q}(G,u)$ is order isomorphic to a rational group. Furthermore, $\quotient{G}{\mathbb{Q}(G,u)}$, as an abstract group, is torsion free. 
\end{proposition}

\proc{Proof.}
The map $\Gamma$ is well-defined. In fact, if $ng=mu$ and $n'g=m'u$, then we get by multiplying with $n'$ and $n$ respectively: $nn'g=n'mu$ and $nn'g=nm'u$. Subtracting we get: $(n'm-nm')u=0$. Since $G$ is torsion free we get
$n'm-nm'=0$, and so $\frac{m}{n} = \frac{m'}{n'}$. Similarly we show that $\Gamma$ is a group homomorphism sending $u$ to 1. If $\Gamma(g)=0$, then $ng=0$ for some $n\in\mathbb{N}$, and so $g=0$ by torsion-freeness of $G$. 
So $\Gamma$ is injective. It is straightforward to show that $\Gamma(g)\geq 0 \Leftrightarrow g\geq 0$, and so $\Gamma$ is an order isomorphism onto its image. Hence $\mathbb{Q}(G,u)$ is order isomorphic to a rational group.

To show that $\quotient{G}{\mathbb{Q}(G,u)}$ is torsion free, let $k\bar{g} = 0$ for some $k\in\mathbb{N}$, where $\bar{g}$ is the image of $g\in G$ under the quotient map. This implies that $kg\in \mathbb{Q}(G,u)$, and so
there exist $n\in\mathbb{N}$, $m\in\mathbb{Z}$ such that $nkg=mu$. Hence $g\in\mathbb{Q}(G,u)$, and so $\bar{g}=0$, thus proving that $\quotient{G}{\mathbb{Q}(G,u)}$ is torsion free. \ep

\proc{Remark.}
Let $(G_1,G_1^+,u_1)\cong (G_2, G_2^+,u_2)$ by a map $\deff{\phi}{G_1}{G_2}$. Then it is easily seen that $\phi\left(\mathbb{Q}(G_1,u_1)\right)= \mathbb{Q}(G_2,u_2)$. So, loosely
speaking, we have that isomorphic dimension groups with distinguished order units have the same rational subgroups. 
\medbreak

The notion of rational subgroup of a dimension group with distinguished order unit depends heavily upon the choice of the order unit as the following example shows. 

\begin{example}
 Let $H$ be a noncyclic subgroup of $\mathbb{Q}$ containing $\mathbb{Z}$. Let $G=H\oplus \mathbb{Z}$ with $G^+ = \set{(h,k)}{h>0, k\in\mathbb{Z}}\cup\{(0,0)\}$. Then $(G,G^+)$ is a simple
dimension group with $\mathrm{Inf}(G)=0\oplus \mathbb{Z}$. If we choose the order unit $u=(1,0)$ for $G$, then one shows easily that $\mathbb{Q}(G,u)=H\oplus0 \cong H$. However, if we choose the 
order unit $\tilde{u}=(1,1)$ for $G$, then $\mathbb{Q}(G,\tilde{u}) = \set{(k,k)}{k\in\mathbb{Z}}\cong \mathbb{Z}$. 
\end{example}

\subsubsection{Ordered Bratteli diagram and the Bratteli-Vershik model} 
\label{subsec:B-V}
An \emph{ordered Bratteli diagram} $(V,E,\geq)$ is a Bratteli diagram $(V,E)$ together with a partial order $\geq$ in $E$ so that edges $e,e'\in E$ are comparable if and only if $r(e)=r(e')$. In other words, 
we have a linear order on each set $r^{-1}(v),\, v\in V\backslash V_0$. We let $E_\text{min}$ and $E_\text{max}$,respectively, denote the minimal and maximal edges if the partially ordered set $E$.

Note that if $(V,E,\geq)$ is an ordered Bratteli diagram and $k<l$ in $\mathbb{Z}^+$, then the set $E_{k+1}\circ E_{k+2}\circ \cdots \circ E_l$ of paths from $V_k$ to $V_l$ with the same range can be given an induced
(lexicographic) order as follows: 
\[(e_{k+1}\circ e_{k+2}\circ \cdots \circ e_l)>(f_{k+1}\circ f_{k+2}\circ \cdots \circ f_l) \]
if for some $i$ with $k+1\leq i\leq l$, $e_j=f_j$ for $i<j\leq l$ and $e_i>f_i$. If $(V',E')$ is a telescoping of $(V,E)$ then, with this induced order from $(V,E,\geq)$, we get again an ordered Bratteli diagram $(V',E',\geq)$.

\begin{definition}
 We say that the ordered Bratteli diagram $(V,E,\geq)$, where $(V,E)$ is a simple Bratteli diagram, is \emph{properly ordered} if
there exists a unique min path $x_\text{min} = (e_1,e_2,\dots)$ and a unique max path $x_\text{max} = (f_1,f_2,\dots)$ in $X_{(V,E)}$. (That is, $e_i\in E_\text{min}$ and $f_i\in E_\text{max}$ for all $i=1,2,\dots$.)
\end{definition}

Let $(V,E)$ be a properly ordered Bratteli diagram, and let $X_{(V,E)}$ be the path space associated to $(V,E)$. Then $X_{(V,E)}$ is a Cantor set. Let $T_{(V,E)}$ be the \emph{lexicographic map} on $X_{(V,E)}$, i.e.\ if 
$x=(e_1,e_2,\dots)\in X_{(V,E)}$ and $x\neq x_\text{max}$ then $T_{(V,E)}x$ is the successor of $x$ in the lexicographic ordering. Specifically, let $k$ be the smallest natural number so that $e_k\notin E_\text{max}$. Let $f_k$ be the 
successor of $e_k$ (and so $r(e_k)=r(f_k)$). Let $(f_1,f_2,\dots,f_{k-1})$ be the unique least element in $E_1\circ E_2 \circ \cdots \circ E_{k-1}$ from $s(f_k)\in V_{k-1}$ to the top vertex $v_0\in V_0$. 
Then $T_{(V,E)}((e_1,e_2,\dots)) = (f_1,f_2,\dots,f_k,e_{k+1},e_{k+2},\dots)$. We define $T_{(V,E)}x_\text{max}=x_\text{min}$. Then it is easy to check that $T_{(V,E)}$ is a minimal homeomorphism on $X_{(V,E)}$. We note that
if $x\neq x_\text{max}$ then $x$ and $T_{(V,E)}x$ are \emph{cofinal}, i.e.\ the edges making up $x$ and $T_{(V,E)}x$, respectively, agree from a certain level on. We will call the Cantor minimal system $(X_{(V,E)},T_{(V,E)})$ a
\emph{Bratteli-Vershik system}. There is an obvious way to telescope a properly ordered Bratteli diagram, getting another properly ordered Bratteli diagram, such that the associated Bratteli-Vershik systems are conjugate -- the map
implementing the conjugacy is the obvious one. By telescoping we may assume without loss of generality that the properly ordered Bratteli diagram has the property that at each level all the minimal edges (respectively the maximal edges)
have the same source.

\begin{theorem}[\cite{HPS}]
\label{th:BVmod}
Let $(X,T)$ be a Cantor minimal system. Then there exists a properly ordered Bratteli diagram $(V,E,\geq)$ such that the associated Bratteli-Vershik system $(X_{(V,E)},T_{(V,E)})$ is conjugate to $(X,T)$. 
\end{theorem}

\proc{Proof sketch.}
Let $x_0\in X$ and let $\{U_n\}_{n\in\mathbb{Z}_+}$ be a decreasing sequence of clopen sets of $X$ such that $U_0=X$ and $U_n\searrow\{x_0\}$. For each $U_n$ we construct a finite number of towers ``built'' over $U_n$. These are determined
by the map $\deff{\lambda_n}{U_n}{\mathbb{N}}$, $\lambda_n(y) = \mathrm{inf}\{m\in\mathbb{N}\, | \, T^my\in U_n\}$. If $\lambda_n(U_n) = \{m_1,m_2,\dots,m_{k_n}\}$, then we get at first $k_n$ towers of height $m_1,m_2,\dots,m_{k_n}$, 
respectively. These may be vertically subdivided, giving rise to more towers (some of them of the same height), such that we obtain the following scenario: The clopen partitions $\{\mathcal{P}_n\}_{n\in\mathbb{Z^+}}$ of $X$ that the towers associated to the various $U_n$'s
generate are nested, $\mathcal{P}_0\prec\mathcal{P}_1\prec\mathcal{P}_2\prec \cdots$, and the union of the $\mathcal{P}_n$'s is a basis for the topology of $X$. We build the properly ordered Bratteli diagram $(V,E, \geq)$ by letting the 
vertices $V_n$ at level $n$ correspond to the various towers built over $U_n$. The ordering of the edges between levels $n-1$ and $n$ is determined by the order in which the towers at level $n$ traverse the towers at level $n-1$.

\proc{Remark.}
The simplest Bratteli-Vershik model $(V,E,\geq)$ for the odometer $(G_\mathfrak{a},T)$ associated to $\mathfrak{a}=(a_i)_{i\in\mathbb{N}}$ is obtained by letting $V_n=1$ for all $n$, and the number of edges between $V_{n-1}$ and
$V_n$ be $a_n$.
\medbreak

\begin{definition}
\label{def:kdyn}
 Let $(X,T)$ be a Cantor minimal system. Define the (additive) \emph{coboundary map} $\deff{\partial_T}{\mathcal{C}(X,\mathbb{Z})}{\mathcal{C}(X,\mathbb{Z})}$ by $\partial_Tf=f-f\circ T^{-1}$. Define 
\[\mathrm{K}^0(X,T)=\quotient{\mathcal{C}(X,\mathbb{Z})}{\partial_T(\mathcal{C}(X,\mathbb{Z}))}\]
and give $\mathrm{K}^0(X,T)$ the induced order, i.e.\ 
\[\mathrm{K}^0(X,T)^+=\set{[f]\in\mathrm{K}^0(X,T)}{[f]=[g], \text{ for some } g\geq 0},\]
where $[f]$ denotes the class of $f\in \mathcal{C}(X,\mathbb{Z})$. Let $\mathds{1} (=[1])$ denote the distinguished order unit corresponding to the constant function $1$ on $X$. We will in the sequel sometimes for short only write
$K^0(X,T)$, the ordering and the order unit being implicitly understood. 
\end{definition}

\proc{Remark.}
Let $(G_\mathfrak{a},\rho_{\hat{\mathds{1}}})$ be the odometer associated to $\mathfrak{a} = (a_i)_{i\in\mathbb{N}}$. Then $K^0(G_\mathfrak{a},\rho_{\hat{\mathds{1}}})$ is order isomorphic to the rational group associated to $\mathfrak{a}$, namely
\[\left\lbrace \left. \frac{m}{a_1a_2\cdots a_n}\right| m\in\mathbb{Z}, n\in\mathbb{N}\right\rbrace,\]
and so, in particular, $\mathbb{Q}(K^0(X,T),\mathds{1}) \cong \mathbb{Q}(K_0(V,E),[\mathbf{1}]$ by a map preserving the canonical order units.
\medbreak

The following result, which is implicit in \cite[Section 2]{GPS}, is highly relevant for this paper and we present a proof due to B. Host (cf.\ \cite[Theorem 2.2]{O}).

\begin{proposition}
 \label{prop:eigenvalue}
 Let $(X,T)$ be a Cantor minimal system and let $p$ be a natural number greater than 1. The following are equivalent:
 \begin{enumerate}[(i)]
  \item\label{pt:1} $\frac{1}{p} \in \mathbb{Q}(K^0(X,T),\mathds{1})$
  \item\label{pt:2} $\frac{2\pi\mathrm{i}}{p}$ is a continuous eigenvalue for $T$, i.e.\ $\exists f\in\mathcal{C}(X)$, $f\neq 0$, such that  $f\circ T = \mathrm{e}^\frac{2\pi\mathrm{i}}{p} f$.
 \end{enumerate}
\end{proposition}

\proc{Proof.} $(\ref{pt:1})\Rightarrow(\ref{pt:2})$: For any set $E$, let $1_E$ denote the characteristic function of $E$. By condition (\ref{pt:1}) there exist continuous functions $\deff{f,g}{X}{\mathbb{Z}}$ such that $pf-1_X=g-g\circ T^{-1}$, where $f$
 is non-zero (and hence $g$ is non-zero). Let $n\in \mathrm{Im}(g)$ and let $A=g^{-1}(\{n\})$. Then $A$ is a non-empty clopen set. For $x\in A$ let $r(x)$ be the first return time of $x$ to $A$, i.e.\ $r(x)$ is the smallest positive
 integer such that $T^{r(x)}x\in A$. Since $T$ is minimal the function $\deff{r}{A}{\mathbb{Z}}$ is uniformly bounded and $\bigcup_{x\in A}\bigcup_{j=1}^{r(x)} T^jx = X$. We get:
 
 \begin{align*}
   p\left(\sum_{j=1}^{r(x)} f(T^jx)\right)-r(x) &= \sum_{j=1}^{r(x)} \left(pf(T^jx)-1_X(T^jx)\right) = \sum_{j=1}^{r(x)} \left(g(T^jx)-g(T^{j-1}x)\right)\\
   & = g(T^{r(x)}x)-g(x).
 \end{align*}

 If $x\in A$, the right hand side is 0, and so the left hand side yields that $r(x)$ must be a (positive) multiple of $p$. Let $\tilde{A}=\bigcup_{l\geq 0} T^{lp}(A)$, and so, in particular, $A$ is an open set. Then for all $x\in A$ the 
 first return time to $\tilde A$ is $p$. Then $X = \sqcup_{j=0}^{p-1} T^j(\tilde A)$ (disjoint union), and so, in particular, $\tilde A, T(\tilde A), \dots, T^{p-1}(\tilde A)$, are clopen sets. Define $f\in\mathcal{C}(X)$ by 
 \[f(x)= \mathrm{e}^\frac{j2\pi\mathrm{i}}{p} \text{ if } x\in T^j(\tilde A).\]
 Then one easily sees that $f\circ T = \mathrm{e}^\frac{2\pi\mathrm{i}}{p}f$, and so $\frac{2\pi\mathrm{i}}{p}$ is a continuous eigenvalue for $T$.
 
$(\ref{pt:2})\Rightarrow(\ref{pt:1})$: Let  $\frac{2\pi\mathrm{i}}{p}$ be a continuous eigenvalue for $T$, and let $f\in\mathcal{C}(X)$ be a non-zero eigenfunction, i.e.\ $f\circ T = \mathrm{e}^\frac{2\pi\mathrm{i}}{p}f$.
 We may assume that $\abs{f(x)}=1$ for all $x\in X$ and that $A= f^{-1}\left(\{1\}\right)$ is non-empty. Now clearly $A, T(A),\dots,T^{p-1}(A)$ are disjoint closed sets, and $T\left(\bigcup_{j=0}^{p-1} T^j (A)\right) = \bigcup_{j=0}^{p-1} T^j (A)$,
 and so by minimality of $T$, we have $\bigcup_{j=0}^{p-1} T^j (A)=X$. Hence $A$ is clopen, and so $1_A\in\mathcal{C}(X)$. Now 
 \begin{align*}
  p1_A-1_X &= (1_A-1_A)+(1_A-1_{T(A)})+\cdots + (1_A-1_{T^{p-1}(A)}\\
  &= 0 + (1_A-1_A\circ T^{-1})+\cdots + (1_A-1_A\circ T^{-(p-1)}).
 \end{align*}
 Also
 \begin{align*}
  &1_A-1_A\circ T^{-j} = \\
  &(1_A-1_A\circ T^{-1}) + (1_A\circ T^{-1}-(1_A\circ T^{-1})\circ T^{-1}) + \cdots + (1_A\circ T^{-(j-1)}-(1_A\circ T^{-(j-1)})\circ T^{-1})
 \end{align*}
 is a coboundary, and so $p1_A-1_X$ is a coboundary. This shows that $\frac{1}{p}\in\mathbb{Q}(K^0(X,T),\mathds{1})$.
$\square$
\medbreak

Combining Proposition \ref{prop:eigenvalue} with Theorems \ref{th:paul} and \ref{th:markley-paul} (and the remarks just preceding these theorems), we get the following result:

\begin{proposition}
\label{prop:maxfactor}
 Let $(X,T)$ be a Cantor minimal system. Then $\mathbb{Q}(K^0(X,T),\mathds{1})$ completely determines the maximal equicontinuous \emph{Cantor} factor of $(X,T)$. This factor is an odometer $(G_\mathfrak{a},\rho_{\hat{\mathds{1}}})$ 
 associated to an $\mathfrak{a}$-adic number $\mathfrak{a} = (a_1, a_2, \dots)$, and $\mathbb{Q}(K^0(X,T),\mathds{1})\cong K^0(G_\mathfrak{a},\rho_{\hat{\mathds{1}}})$, where the latter group is described in the Remark preceding 
 \ref{prop:eigenvalue}. In particular, if $(X,T)$ is a Toeplitz flow then $\mathbb{Q}(K^0(X,T),\mathds{1})$ completely determines the maximal equicontinuous factor of $(X,T)$.
\end{proposition}

The following theorem is a fairly straightforward corollary of Theorem \ref{th:BVmod}.

\begin{theorem}[\cite{HPS}]
\label{th:realization}
 Let $(X,T)$ be a Cantor minimal system. Then $(\mathrm{K}^0(X,T), \mathrm{K}^0(X,T)^+)$ is a simple dimension group, and all (noncyclic) simple dimension groups $(G,G^+,u)$ with distinguished order unit $u$ are order isomorphic to $(\mathrm{K}^0(X,T), \mathrm{K}^0(X,T)^+,\mathds{1})$
for some Cantor minimal system $(X,T)$. In particular, if $(V,E,\geq)$ is a properly ordered Bratteli diagram such that $(X,T)$ is conjugate to $(X_{(V,E)}, T_{(V,E)})$, then $K^0(X,T)\cong K_0(V,E)$ (as ordered groups with 
canonical order units), and so, in particular, $\mathbb{Q}(K^0(X,T),\mathds{1}) \cong \mathbb{Q}(K_0(V,E),[\mathbf{1}])$. Furthermore, the state space $S_\mathds{1}(K^0(X,T))$ may be identified in an obvious way with the Choquet simplex $M(X,T)$ of $T$-invariant probability measures; in fact, these two Choquet simplices are affinely 
homeomorphic.
\end{theorem}

\proc{Remark.}
Changing the order unit corresponds dynamically to considering induced systems. (Cf.\ Definitions \ref{def:induced} and \ref{def:kaku}.) In fact, let $(V,E,\geq)$ be a simply ordered Bratteli diagram, and let $(V',E',\geq)$ be the resulting simply ordered Bratteli diagram after we have
made a finite change to $(V,E,\geq)$. (By a finite change we mean adding and/or removing a finite number of edges and then making arbitrary choices of linear orderings of the edges meeting at the same vertex for a finite number of
vertices.) Then $(X_{(V,E)}, T_{(V,E)})$ is Kakutani equivalent to $(X_{(V',E')}, T_{(V',E')})$, which can be seen as an immediate consequence of how the Vershik map is defined. Furthermore, if $(X,T)$ is Cantor minimal with associated
dimension group $K^0(X,T)$ with distinguished order unit, then choosing a new order unit, say $u$, there exists a Cantor minimal system $(Y,S)$ which is Kakutani equivalent to $(X,T)$ such that 
$(K^0(Y,S),K^0(Y,S)^+,\mathds{1})\cong (K^0(X,T),K^0(X,T)^+,u)$. In fact, $(Y,S)$ is obtained from $(X,T)$ by making a finite change to the Bratteli-Vershik model for $(X,T)$.
\medbreak

\begin{theorem}[\cite{GPS}]
\label{th:kgroups}
 The Cantor minimal systems $(X,T)$ and $(Y,S)$ are \emph{strong orbit equivalent} if and only if $(K^0(X,T),K^0(X,T)^+,\mathds{1})\cong (K^0(Y,S),K^0(Y,S)^+,\mathds{1})$.
\end{theorem}

\proc{Remark.}
The idea behind introducing $K^0(X,T)$ with an ordering for a Cantor minimal system $(X,T)$ comes from (non-commutative) operator algebra theory. In fact, one can show that $K^0(X,T)$ as defined above is abstractly isomorphic to
the $K_0$-group of the associated $C^\ast$-crossed product $\cont{X}\rtimes_T \mathbb{Z}$. The latter $K_0$-group comes with a natural order, which translates to the order we introduced on $K^0(X,T)$.
It turns out that the ordered $K_0$-group of $\cont{X}\rtimes_T \mathbb{Z}$, with its natural scaling corresponding to the unit element, is a complete isomorphism invariant, and so it follows that $(X,T)$ is strong orbit
equivalent to $(Y,S)$ if and only if $\cont{X}\rtimes_{T} \mathbb{Z}$ is $\ast$-isomorphic to $\cont{Y}\rtimes_{S} \mathbb{Z}$ \cite{GPS}.
\medbreak

\begin{theorem}[\cite{GPS}]
\label{th:orbit}
The Cantor minimal systems $(X,T)$ and $(Y,S)$ are \emph{orbit equivalent} if and only if $(\tilde{K^0(X,T)},\tilde{K^0(X,T)^+},\tilde{\mathds{1}})\cong (\tilde{K^0(Y,S)},\tilde{K^0(Y,S)^+},\tilde{\mathds{1}})$, where 
\newline $\tilde{K^0(X,T)} = \quotient{K^0(X,T)}{\mathrm{Inf}(K^0(X,T))}$, $\tilde{K^0(Y,S)} = \quotient{K^0(Y,S)}{\mathrm{Inf}(K^0(Y,S))}$, and the positive cones and order units are induced by the quotient maps.
\end{theorem}

Let $(V,E,\geq)$ be a simply ordered Bratteli diagram, and let $(X_{(V,E)},T_{(V,E)})$ be the associated Bratteli-Vershik system. For each $k\in \mathbb{N}$ let $P_k$ as above denote the paths from $V_0$ to $V_k$, i.e.\ the paths from $v_0\in V_0$ 
to some $v\in V_k$. For $x\in X_{(V,E)}$ we associate the bi-infinite sequence $\pi_k(x) = \left(\tau_k(T_{(V,E)}^n x)\right)_{n=-\infty}^\infty\in P_k^\mathbb{Z}$ over the finite alphabet $P_k$, where $\deff{\tau_k}{X_{(V,E)}}{P_k}$ 
is the truncation map. Let $S_k$ denote the shift map on $P_k^\mathbb{Z}$. Then the following diagram commutes
\begin{center}
 \includegraphics{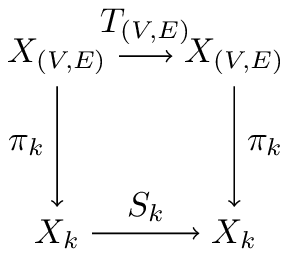}
\end{center}
where $X_k = \pi_k(X_{(V,E)})$. One observes that $\pi_k$ is a continuous map, and so $X_k$ is a compact shift-invariant subset of $P_k^\mathbb{Z}$. So $(X_k,S_k)$ is a factor of $(X_{(V,E)}, T_{(V,E)})$. For $k>l$ there is an 
obvious factor map $\deff{\pi_{k,l}}{X_k}{X_l}$, and one can show that $(X_{(V,E)}, T_{(V,E)})$ is the inverse limit of the system $\{(X_k,S_k)\}_{k\in\mathbb{N}}$. All the systems $(X_k,S_k)$ are clearly expansive. One has
the following result which will be important for us.

\begin{proposition}
\label{prop:expansive-conj}
Assume $(X_{(V,E)}, T_{(V,E)})$ is expansive. Then there exists $k_0\in \mathbb{N}$ such that for all $k\geq k_0$, $(X_{(V,E)}, T_{(V,E)})$ is conjugate to $(X_k,S_k)$ by the map $\deff{\pi_k}{X_{(V,E)}}{X_k}$.
\end{proposition}

\proc{Proof.}
Since the $\pi_k$'s are factor maps, all we need to show is that there exists $k_0$ such that $\pi_k$ is injective for all $k\geq k_0$. Recall that $(X_{(V,E)}, T_{(V,E)})$ being expansive means that there exists $\delta>0$ such that given 
$x\neq y$ there exists $n_0\in \mathbb{Z}$ such that $d(T_{(V,E)}^{n_0}x,T_{(V,E)}^{n_0}y)>\delta$, where $d$ is some metric on $X_{(V,E)}$ compatible with the topology. Choose $k_0$ such that $d(x,y)<\delta$ if $x$ and $y$ agree (at least) on the $k_0$ first edges. Now assume that 
$\pi_k(x)=\pi_k(y)$ for some $k\geq k_0$. By the definition of $\pi_k$ this means that, for all $n\in \mathbb{Z}$, $\tau_k(T_{(V,E)}^n x) = \tau_k(T_{(V,E)}^n y)$, and so $d(T_{(V,E)}^n x,T_{(V,E)}^n y)<\delta$ for all $n\in\mathbb{Z}$ because of our choice of $k_0$. 
This contradicts that $d(T_{(V,E)}^{n_0}x,T_{(V,E)}^{n_0}y)>\delta$. Hence $\pi_k$ is injective for all $k\geq k_0$, proving the theorem.
$\square$
\medbreak

\subsection{Order embeddings of dimension groups associated to factor maps.}

We have seen that  Toeplitz flows can be characterized by being (expansive) almost one-to-one Cantor minimal extensions of odometer systems (cf.\ Theorem \ref{th:markley-paul}). The following theorem
will therefore be important for us since it relates the $K^0$-groups of the extension and the factor, respectively, of Cantor minimal systems. 

\begin{theorem}[\protect{\cite[Proposition 3.1]{GW}}]
\label{th:GW}
Let $(X,T)$ and $(Y,S)$ be Cantor minimal systems such that $(X,T)$ is an extension of $(Y,S)$ by the factor map $\deff{\pi}{X}{Y}$. Then $\deff{\pi^\ast}{K^0(Y,S)}{K^0(X,T)}$ defined by $\pi^\ast([h]) = [h\circ \pi]$
is an order embedding, i.e.\ $[h]\geq 0$ if and only if $\pi^\ast([h])\geq 0$ for $h\in \mathcal{C}(Y,\mathbb{Z})$. (Here $[h]$ and $[h\circ\pi]$ denote the class of $h$ and $h\circ \pi$ in $K^0(Y,S)$ and $K^0(X,T)$, respectively.)
\end{theorem}

\proc{Remark.}
The proof of Theorem \ref{th:GW} has two ingredients. The first is the use of the Gottschalk-Hedlund lemma \cite{GH}, which in our context says that $g\in \mathcal{C}(X,\mathbb{Z})$ is a coboundary, i.e.\ $g=f-f\circ T^{-1}$
for some $f\in\mathcal{C}(X,\mathbb{Z})$, if and only if $\sup_n \left| \sum_{i=0}^{n-1} g(T^ix_0)\right| < \infty $ for some $x_0\in X$. This will establish that $\pi^\ast$ is well-defined. 
The second ingredient, which gives the order embedding, is applying Theorem \ref{th:simplex} (cf.\ also Theorem \ref{th:realization}) together with the fact that $\pi$ induces a surjective map of $M(X,T)$ onto $M(Y,S)$.

\section{Proofs of the main results}
The proofs of Theorem \ref{th:twoequalsets}, Theorem \ref{th:realizing} and Theorem \ref{th:threeequalsets} stated in Section \ref{sec:results} will rest heavily upon results obtained earlier by Gjerde-Johansen \cite{gjerde-johansen} and Sugisaki \cite{sugisaki-3}, \cite{sugisaki-1}, \cite{sugisaki-2}, where the paper \cite{sugisaki-2} is extending and being inspired, so to speak, by an analogous result proved by Giordano, Putnam, Skau in \cite{GPS2}. The proofs in \cite{sugisaki-3}, \cite{sugisaki-1} and \cite{sugisaki-2}
are rather technical, involving very clever manipulations of Bratteli diagrams. Looking carefully at crucial steps in the proofs of the main theorems in \cite{sugisaki-1} and \cite{sugisaki-2}, in particular, we could deduce more
specific properties of the Bratteli diagrams that appear, starting with our basic setting. This, combined with the Bratteli-Vershik model for Toeplitz flows established in \cite{gjerde-johansen} and the embedding
result of \cite{GW}, will, loosely speaking, give the proofs of our three first theorems. (The proof of the fourth theorem, Theorem \ref{th:real_simplex}, requires a somewhat different approach.) However, this should not be construed as saying that our theorems are simply corollaries of these earlier results. But one could perhaps say that 
the results we obtain represent the culmination of the study of (topological) Toeplitz flows from the perspective of orbit equivalence and/or $K^0$-groups.

\begin{lemma}
\label{lem:ers}
 The Bratteli--Vershik model associated to a properly ordered Bratteli diagram $(V,E, \geq)$ with the ERS-property has the odometer corresponding to the supernatural number of $(V,E)$  
 (cf.\ Definition \ref{def:ers}) as a factor. The factor map $\pi$ is almost one-to-one. If the Bratteli-Vershik system $(X_{(V,E)},T_{(V,E)})$ is expansive, then for some $k$ we get that $\pi_k(x_{\mathrm{min}})$ 
 is a Toeplitz sequence and $\deff{\pi_k}{X_{(V,E)}}{X_k}$ is a conjugacy map between $(X_{(V,E)}, T_{(V,E)})$ and $(X_k,S_k)$ (cf.\ Proposition \ref{prop:expansive-conj}). In particular, $(X_{(V,E)},T_{(V,E)})$ is a Toeplitz flow.
\end{lemma}

\proc{Proof.}

\begin{figure}
\centering 
\includegraphics[scale=0.8]{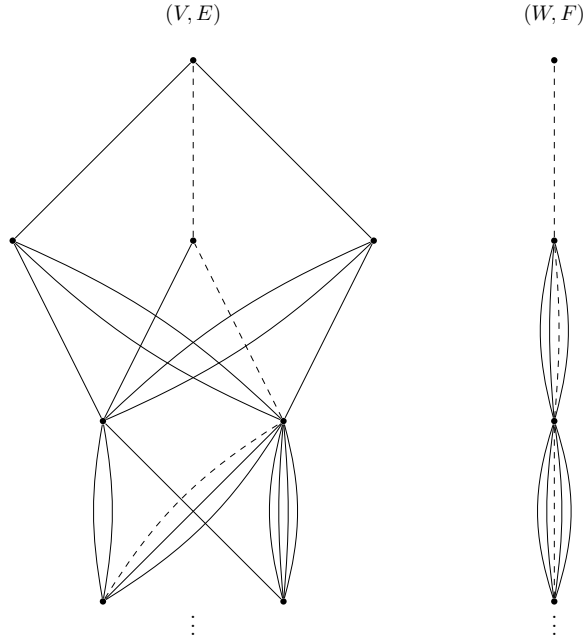}
\caption{Assuming the edges in both $(V,E)$ and $(W,F)$ are ordered left to right, the factor map $\pi$ will take the dashed path in $(V,E)$ to the dashed path in $(W,F)$.}
 \label{fig:ersfactor}
\end{figure}
We will give the proof by appealing to Figure \ref{fig:ersfactor}, where we have drawn on the left a particular Bratteli diagram $(V,E)$ with the ERS-property. Let the incidence matrices be $(M_n)_{n=1}^\infty$, say, 
with $r_n$ equal to the (constant) row sum of $M_n$. We ``collapse'' $(V,E)$ in an obvious way to create the diagram $(W,F)$ on the right in Figure \ref{fig:ersfactor}, where $|W_n|=1$ for all $n$ and the number of 
edges between $W_{n-1}$ and $W_n$ is $r_n$. Let them be linearly ordered as $f_1^{(n)}<f_2^{(n)}<\cdots<f_{r_n}^{(n)}$. The factor map $\deff{\pi}{X_{(V,E)}}{X_{(W,F)}}$ is defined by reading off the labels, so to say.
That is, let $x=(e_n)_{n=1}^\infty\in X_{(V,E)}$ and let $e_n$ be the $k_n$'th edge in the linear ordering of $r^{-1}(e_n)$. We then map $e_n$ to $f_{k_n}^{(n)}$, thus getting $\pi(x)\in  X_{(W,F)}$. It is easy to see that 
$\pi$ is a continuous map that intertwines $T_{(V,E)}$ and $T_{(W,F)}$. Also, $\pi^{-1}\{\pi(x_\text{min})\} = \{x_\text{min}\}$, where $x_\text{min}$ is the unique minimal path of $X_{(V,E)}$, and so $\pi$ is almost one-to-one.
The map $\pi$ is onto since the image of $\pi$ is a compact, hence closed, subset of $X_{(W,F)}$, and $\pi$ clearly maps $\text{orbit}_{T_{(V,E)}}(x_\text{min})$ to a dense subset of $X_{(W,F)}$. The factor $(X_{(W,F)}, T_{(W,F)})$
is by construction an odometer with the properties stated. \linebreak
If $(X_{(V,E)},T_{(V,E)})$ is expansive, then for some $k\in\mathbb{N}$, $\deff{\pi_k}{X_{(V,E)}}{X_k}$ is a conjugacy map by Proposition \ref{prop:expansive-conj}. As pointed out in Section \ref{subsec:B-V} we may assume that at each level
all the minimal edges (respectively the maximal edges) have the same source. This scenario is obtained by telescoping the original diagram, and all essential properties are preserved by this operation.
The ERS-property of $(V,E)$ implies that $\pi_k(x_{\mathrm{min}})$ is a Toeplitz sequence, a fact that is easily shown; we omit the details. So $(X_{(V,E)},T_{(V,E)})$ is a Toeplitz flow 
(which, incidentally, also is a consequence of Theorem \ref{th:markley-paul}).
\ep\medbreak

\proc{Remark.}
Lemma \ref{lem:ers} is the easy part of Theorem 8 of \cite{gjerde-johansen}. Conversely, if one starts with a Toeplitz flow $(X,T)$ one can construct a properly ordered Bratteli diagram $(V,E,\geq)$ such that $(V,E)$ has the ERS-property
and $(X,T)\cong (X_{(V,E)},T_{(V,E)})$. This is achieved by merging the structure of Toeplitz flows as described in Section \ref{sec:toeplitz} with the construction described in the proof sketch of Theorem \ref{th:BVmod}. 
(See \cite[Theorem 8]{gjerde-johansen}.) In \cite{sugisaki-3} it is proved that for any ordered Bratteli diagram $(V,E,\geq)$ such that $(V,E)$ has the ERS-property, there exists another properly ordered Bratteli diagram $(V',E',\geq)$
such that $(V',E')$ has the ERS-property, with $K_0(V',E')\cong K_0(V,E)$ (as ordered groups with distinguished order units) and $(X_{(V',E')},T_{(V',E')})$ is Toeplitz. By Theorem \ref{th:realization}, $K^0(X,T)\cong K^0(X_{(V',E')},T_{(V',E')})$, 
and so $(X,T)$ is strong orbit equivalent to $(X_{(V',E')},T_{(V',E')})$ by Theorem \ref{th:kgroups}.
\medbreak

\begin{lemma}
\label{lem:torsionfree}
 Let $\deff{\pi}{(X,T)}{(Y,S)}$ be an almost one-to-one factor map between Cantor minimal systems. Then
\[\quotient{K^0(X,T)}{\pi^\ast(K^0(Y,S))}\]
is torsion free. 
\end{lemma}

\proc{Proof.}
We will give a heuristic argument in order to highlight the basic idea behind the proof, avoiding the somewhat messy details that tend to obscure the understanding. (For a more detailed proof, see \cite[Theorem 3.1]{sugisaki-2}.)
Since $\pi$ is a factor map, we know from Theorem \ref{th:GW} that $\deff{\pi^\ast}{K^0(Y,S)}{K^0(X,T)}$ is an order embedding where $\pi^\ast([h]) = [h\circ \pi]$. (Here $h\in\mathcal{C}(Y,\mathbb{Z})$ and $[h]$ and $[h\circ\pi]$ 
denote the classes of $h$ and $h\circ \pi$ in $K^0(Y,S)$ and $K^0(X,T)$, respectively.)
We will construct specific Bratteli-Vershik models for $(X,T)$ and $(Y,S)$, respectively. Let $x_0\in X$ be such that $\pi^{-1}(\pi(x_0)) = x_0$. Let $\{U_n\}_{n\in\mathbb{Z}^+}$ be a decreasing sequence of clopen sets in $Y$
such that $U_0=Y$ and $U_n\searrow\{y_0\}$, where $y_0=\pi(x_0)$. We now proceed as described in the proof sketch of Theorem \ref{th:BVmod}. Let $\{\mathcal{P}_n\}_{n\in\mathbb{Z}^+}$ be the nested sequence of clopen partitions of $Y$ 
 associated to the tower constructions built over the various $U_n$'s, such that the union of the $\mathcal{P}_n$'s is a basis for the topology of $Y$. Let $(W,F,\geq)$ be the resulting properly ordered Bratteli diagram, from which we get
 a Bratteli-Vershik model for $(Y,S)$. Consider the clopen sets $\tilde{U_n} = \pi^{-1}(U_n)$, $n\in\mathbb{Z}^+$. Clearly $\{\tilde{U_n}\}_{n\in\mathbb{Z}^+}$ is a descending sequence of clopen subsets of $X$ such that $\tilde{U_0}=X$,
 $\tilde{U_n}\searrow\{x_0\}$. Proceeding again as described in the proof sketch of Theorem \ref{th:BVmod} we get a simply ordered Bratteli diagram $(V, E, \geq)$, from which we get a Bratteli-Vershik model for $(X,T)$.
 Now $K^0(X,T)\cong K_0(V,E)$ and $K^0(Y,S)\cong K_0(W,F)$ as ordered groups with canonical order units. We note that the functions $\deff{\lambda_n}{U_n}{\mathbb{N}}$ and $\deff{\tilde \lambda_n}{\tilde U_n}{\mathbb{N}}$ defined by
 $\lambda_n(u) = \mathrm{inf}\{m\in\mathbb{N}\, | \, S^m u \in U_n\}$ and $\tilde\lambda_n(\tilde u) = \mathrm{inf}\{m\in\mathbb{N}\, | \, T^m \tilde u \in \tilde U_n\}$, respectively, are related by $\lambda_n(u) = \tilde \lambda_n (\tilde u)$
 if $\tilde u \in \pi^{-1}(u)$. This has implications for how $K_0(W,F)$ is embedded in $K_0(V,E)$. Loosely speaking, each $w\in W_n$ (which corresponds to a tower $T_w$ over $U_n$) is split into a finite number of vertices $v_1, v_2,\dots, v_l$
 in $V_n$ (this corresponds to the tower $T_w$ being subdivided into $l$ towers $T_{v_1},T_{v_2},\dots,T_{v_l}$ over $\tilde U_n$, each of the same height as $T_w$). The factor map $\deff{\pi}{X_{(V,E)}}{X_{(W,F)}}$ is a kind of 
 ``collapsing'' map similarly to the one exhibited in Figure \ref{fig:ersfactor}. While the scenario exhibited in Figure \ref{fig:ersfactor} is very simple, it does illustrate the essential point. In fact, the single vertices at levels 1, 2
 and 3 of $(W,F)$ in Figure \ref{fig:ersfactor} split into three, two and two vertices, respectively, at levels 1, 2 and 3 of $(V,E)$. The image of the group element of $K_0(W,F)$ that is $-7$, say, at level 2 of $(W,F)$ in Figure \ref{fig:ersfactor}
 is represented by $(-7,-7)$ at level 2 of $(V,E)$. In general, a group element in $K_0(W,F)$ which is represented as being $b\in\mathbb{Z}$ at $w\in W_n$, and zero at the other vertices in $W_n$, is mapped to the group element in 
 $K_0(V,E)$ that is represented by being $b$ at each of the vertices in $V_n$ associated to $w$, and zero elsewhere. This extends by linearity in an obvious way to any element in $K_0(W,F)$ that are represented as a vector at level $n$. This
 ``locally constancy'' property, which is preserved at the higher levels under the canonical mappings of the Bratteli diagram, is what completely characterizes the embedding of $K_0(W,F)$ into $K_0(V,E)$. This clearly implies that
 $\quotient{K_0(V,E)}{\pi^\ast(K_0(W,F))}$ is torsion free. If namely $g\in K_0(V,E)$ is such that $kg\in \pi^\ast(K_0(W,F))$ for some $k\in\mathbb{N}$, then $kg$ is represented by some ``locally constant'' vector at level $n$ of $(V,E)$,
 say. But then clearly $g$ is also represented by a ``locally constant'' vector (at some higher level of $(V,E)$ than $n$, perhaps), and hence $g$ lies in $\pi^\ast(K_0(W,F))$.
\ep
\medbreak

\proc{Remark.}
The converse of Lemma \ref{lem:torsionfree} is not true. In a recent paper by Glasner and Host they construct Cantor minimal systems $(Y,S)$ and $(X,T)$ such that $(X,T)$ is an extension of $(Y,S)$ by a map
$\deff{\pi}{(X,T)}{(Y,S)}$ which is not an almost one-to-one extension, and $\quotient{K^0(X,T)}{\pi^\ast(K^0(Y,S))}$ is non-zero and torsion free \cite[Appendix C]{glasner-host}. In fact, their example can be adjusted to make $(Y,S)$
an odometer and $(X,T)$ to be expansive.
\medbreak

\begin{lemma}
 \label{lem:super}
 Let $(V,E, \geq)$ be a simple Bratteli diagram with the ERS-property. Let $M_n$ be the incidence matrix between levels $n-1$ and $n$, and let $r_n$ be the (constant) row sum of $M_n$. Let $m$ be the supernatural number $\prod_{n=1}^\infty r_n$. Let 
 $H = \mathbb{Q}(K_0(V,E),[{\bf 1}])$ where $[{\bf 1}]$ is the canonical order unit of $K_0(V,E)$. Then $H$ is order isomorphic to the rational group $G(m)$ by a map sending $[{\bf 1}]$ to $1\in G(m)$. Furthermore, $H$ is represented in an obvious way by constant vectors at each level of $(V,E)$, i.e.\ vectors of the form $(a,a,a,\dots,a)^{\mathrm{tr}}\in \mathbb{Z}^{|V_n|}$ for each $n$ ($\mathrm{tr}$ denotes the transpose).
\end{lemma}

\proc{Proof.}
The proof is an immediate consequence of the fact that 
\[M_n(1,1,\dots,1)^{\mathrm{tr}} = (r_n,r_n,\dots,r_n)^{\mathrm{tr}}\]
where $(1,1,\dots,1)^{\mathrm{tr}}\in\mathbb{Z}^{|V_{n-1}|}$, $(r_n,r_n,\dots,r_n)^{\mathrm{tr}}\in\mathbb{Z}^{|V_{n}|}$, which yields
\[M_nM_{n-1}\cdots M_1(1) = \left(\prod_{k=1}^{n}r_k,\prod_{k=1}^{n}r_k,\dots,\prod_{k=1}^{n}r_k\right)^{\mathrm{tr}}\in\mathbb{Z}^{|V_{n}|}.\]
\ep\medbreak

\proc{Remark.}
We will say that $(V,E)$ is an ERS realization of $G\cong K_0(V,E)$ with respect to a subdimension group $H\subseteq G$ if $H$ is embedded in $K_0(V,E)$ as in Lemma \ref{lem:super}.
\medbreak

\begin{lemma}
 \label{lem:sugisaki-extension}
 Let $(J,J^+,1)$ be a noncyclic rational group (cf.\ Definition \ref{def:rational}), and let $(Y,S)$ be an odometer such that $K^0(Y,S)\cong J$ (as ordered groups with distinguished order units). Let $\deff{\iota}{J}{G}$ be an order
 embedding of $J$ into a simple dimension group $(G,G^+,u)$ preserving the order units, such that $\quotient{G}{\iota(J)}$ is torsion free. There exists a properly ordered Bratteli diagram $(V,E,\geq)$ such that:
 \begin{enumerate}[(i)]
  \item \label{ext:1} $(V,E)$ has the ERS-property
  \item \label{ext:2} $K^0(X_{(V,E)},T_{(V,E)}) \cong G$ (as ordered groups with distinguished order units).
  \item \label{ext:3} $(X_{(V,E)},T_{(V,E)})$ is an almost one-to-one extension of $(Y,S)$.
  \item \label{ext:4} $(Y,S)$ is the maximal equicontinuous factor of $(X,T)$.
  \item \label{ext:5} $(J,J^+,1) \cong \mathbb{Q}(K_0(V,E),[{\bf 1}])$.
 \end{enumerate}
\end{lemma}

\proc{Proof.}
Assertions (\ref{ext:2}) and (\ref{ext:3}) are the main result of \cite{sugisaki-2}, namely Theorem 1.1 (see also Corollary 1.2). Since the properly ordered Bratteli diagram associated to $(Y,S)$ is very special --
having a single vertex at each level -- the almost one-to-one extension of $(Y,S)$ constructed in the proof of Theorem 1.1 in \cite{sugisaki-2}, which is obtained by constructing a properly ordered Bratteli diagram $(V,E,\geq)$ based on the one
associated to $(Y,S)$, will have property (\ref{ext:1}). (For details, cf.\ Remark 3.2 and Proposition 3.3 of \cite{sugisaki-2}.) The assertion (\ref{ext:4}) is a consequence of Theorem \ref{th:paul}. The assertion (\ref{ext:5}) follows from (\ref{ext:3}) and Lemma \ref{lem:super}. In fact, $K^0(Y,S)$ embeds into $K_0(V,E)$ as constant vectors at each level of $(V,E)$. (Cf.\ the proof of Lemma \ref{lem:torsionfree}, keeping in mind that in our case the properly ordered Bratteli diagram $(W,F,\geq)$ appearing there and being associated to $(Y,S)$, has one vertex at each level.)
\ep
\medbreak

\begin{lemma}
 \label{lem:sugisaki-expansive}
 Let $(V,E,\geq)$ be a properly ordered Bratteli diagram such that $(V,E)$ has the ERS-property. Let $0\leq t < \infty$. There exists a properly ordered Bratteli diagram $(\tilde V,\tilde E,\geq)$ such that:
 \begin{enumerate}[(i)]
  \item \label{ext:1} $(\tilde V,\tilde E)$ has the ERS-property
  \item \label{ext:2} $K_0(\tilde V,\tilde E)\cong K_0(V,E)$ (as ordered groups with distinguished order units).
  \item \label{ext:3} $(X_{(\tilde V,\tilde E)},T_{(\tilde V,\tilde E)})$ is expansive.
  \item \label{ext:4} $h(T_{(\tilde V,\tilde E)})=t$.
 \end{enumerate}
\end{lemma}

\proc{Proof.}
The assertions (\ref{ext:2}), (\ref{ext:3}) and (\ref{ext:4}) are the main result (Theorem 1.1) of \cite{sugisaki-1}. (Note that by Theorems \ref{th:realization}, \ref{th:kgroups} and Proposition \ref{prop:subshift}, respectively,
strong orbit equivalence is related to $K_0$-groups, and expansiveness is related to subshifts, respectively.) Now in the proof of Theorem 1.1 of \cite{sugisaki-1} various simple Bratteli diagrams are constructed, modifying the given
Bratteli diagram $(V,E)$. However, each modification preserves the ERS-property of the original Bratteli diagram $(V,E)$. (For details, cf.\ Propositions 4.2, 4.4 and Sections 5.1, 5.3 and 5.4 of \cite{sugisaki-1}.) So the properly ordered 
Bratteli diagram $(\tilde V, \tilde E,\geq)$ that eventually arises in the proof of Theorem 1.1 of \cite{sugisaki-1} will have all the properties listed in Lemma \ref{lem:sugisaki-expansive}.
\ep\medbreak

\proc{Proof of Theorem \ref{th:twoequalsets}.}
Let $(X,T)$ be a Toeplitz flow (with $h(T)=t$). By Theorem \ref{th:markley-paul}, $(X,T)$ is an almost one-to-one extension of an odometer $(Y,S)$, the factor map being $\deff{\pi}{X}{Y}$. By Theorem \ref{th:GW}, 
$\deff{\pi^\ast}{K^0(Y,S)}{K^0(X,T)}$ is an order embedding sending the distinguished order unit of $K^0(Y,S)$ to the one in $K^0(X,T)$. Set $G=K^0(X,T)$, $H=\pi^\ast(K^0(Y,S))$. Then $H$ is a noncyclic subgroup of $G$ of rank one such that
$H\cap G^+ \neq \{0\}$.

Conversely, assume that $(G,G^+)$ is a simple dimension group containing the noncyclic subgroup $H$ such that $H\cap G^+ \neq \{0\}$. Let $u\in H\cap G^+$ be any non-zero element. We consider the simple dimension group $(G,G^+,u)$ with
distinguished order unit $u$. Now $H$ is a subgroup of the rational subgroup $\mathbb{Q}(G,u)$ of $G$. In fact, if $h\in H$ then there exists $m,n\in\mathbb{Z}$ such that $nh=mu$, since $H$ is of rank one. (We may assume without loss
 of generality that $n\in\mathbb{N}$.) In particular, we get that $\mathbb{Q}(G,u)$ is noncyclic. Now we apply Lemma \ref{lem:sugisaki-extension} with $J=\mathbb{Q}(G,u)$ and $\deff{\iota}{J}{G}$ the inclusion map. We keep the notation
 of Lemma \ref{lem:sugisaki-extension}. Apply Lemma \ref{lem:sugisaki-expansive} to the properly ordered Bratteli diagram $(V,E,\geq)$ constructed in Lemma \ref{lem:sugisaki-extension} to get $(\tilde V, \tilde E,\geq)$ with the properties listed in Lemma \ref{lem:sugisaki-expansive}. Set $X= X_{(\tilde V, \tilde E)}$ and $T=T_{(\tilde V, \tilde E)}$. 
 Then $K^0(X,T)\cong K_0(\tilde V, \tilde E)\cong K_0(V,E)\cong G$ as ordered groups with distinguished order units. This completes the proof of Theorem \ref{th:twoequalsets}.
\ep\medbreak

\proc{Proof of Theorem \ref{th:realizing}.}
Set $J=\mathbb{Q}(G,u)$. By assumption $J$ is a noncyclic rational group (hence of rank one). The inclusion map $\deff{\iota}{J}{G}$ is an order embedding, and by Proposition \ref{prop:quotient}, $\quotient{G}{J}$ is torsion free.
Let $(Y,S)$ be an odometer system such that $K^0(Y,S)\cong J$ (as ordered groups with distinguished order units). Applying Lemma \ref{lem:sugisaki-extension} we get a  properly ordered Bratteli diagram $(V,E,\geq)$ satisfying the properties
listed in Lemma \ref{lem:sugisaki-extension}. Now apply Lemma \ref{lem:sugisaki-expansive} to $(V,E,\geq)$ to get $(\tilde V, \tilde E,\geq)$ satisfying the properties listed in Lemma \ref{lem:sugisaki-expansive}.
Let $(X,T)=(X_{(\tilde V, \tilde E)}, T_{(\tilde V, \tilde E)})$. By Lemma \ref{lem:ers} we get that $(X,T)$ is a Toeplitz flow. Finally, invoking Proposition \ref{prop:eigenvalue}, we get that $(X,T)$ satisfies all the properties listed in Theorem \ref{th:realizing}. (Recall that by Remark to Proposition \ref{prop:quotient}
we have $(\mathbb{Q}(G,u)\cong)\, \mathbb{Q}(K_0(V,E),[\mathbf{1}]) \cong \mathbb{Q}(K_0(\tilde V,\tilde E),[\mathbf{1}])$.)
\ep
\medbreak

\proc{Proof of Theorem \ref{th:threeequalsets}.}
By Theorem \ref{th:realizing} we get the inclusion $\mathcal{G}\subseteq \mathcal{T}_t$. The inclusion $\mathcal{T}_t\subseteq \mathcal{G}$ is a consequence of Theorem \ref{th:markley-paul} and Proposition \ref{prop:maxfactor}. 
Hence we have $\mathcal{G} = \mathcal{T}_t$. 

The inclusion $\mathcal{B}\subseteq \mathcal{G}$ follows from Lemma \ref{lem:super}. The inclusion $\mathcal{T}_t\subseteq \mathcal{B}$ follows from the main result, Theorem 8, of \cite{gjerde-johansen}. Altogether this implies
that $\mathcal{T}_t=\mathcal{G}= \mathcal{B}$.
\ep
\medbreak

\proc{Proof of Theorem \ref{th:real_simplex}.}
We prove the first part: Let $H$ be a dense rational subgroup of $\mathbb{Q}$. Let $x_0$ be any point in $K$, and let $\mathcal{E}$ be a countable set in $\mathrm{Aff}(K)$, the continuous, real-valued affine functions on $K$, such that $\mathcal{E}$ is a dense
(in the uniform topology) subset of $\set{a\in \mathrm{Aff}(K)}{a(x_0) = 0}$. Let $G$ be the (countable) additive subgroup of $\mathrm{Aff}(K)$ generated by $H$ and $\mathcal{E}$, where we identify every element in $H$ with a constant
affine function. Then $(G,G^+,1)$ is a simple dimension group with order unit the constant function 1, and $G^+$ is the obvious positive cone; furthermore, $\mathrm{Inf}(G)=\{0\}$, and the state space $S_1(G)$ is affinely
homeomorphic with $K$ (cf.\ Theorem \ref{th:simplex}). It is a simple matter to check that $\mathbb{Q}(G,1)=H$. By Theorem \ref{th:realizing} there exists a Toeplitz flow $(X,T)$ such that $(G,G^+,1)\cong (K^0(X,T),K^0(X,T)^+,\mathds{1})$ and $h(T)=t$. By Theorem \ref{th:realization}
we get that $M(X,T)$ is affinely homeomorphic to $S_1(G)$, and so to $K$. Now if $(G_1, G_1^+, 1)\cong (G_2, G_2^+, 1)$, where $G_1$ and $G_2$ are constructed as $G$ above, then $\mathbb{Q}(G_1,1)\cong\mathbb{Q}(G_2,1)$ as ordered groups by a positive map sending 1 to 1. Clearly one can choose an uncountable
family of non-isomorphic groups $H$ of the type described above. Hence there exists an uncountable family of non-isomorphic dimension groups $(G,G^+,1)$ of the type constructed above. The corresponding uncountable family
of Toeplitz flows $(X,T)$ are then pairwise non-orbit equivalent by Theorem \ref{th:orbit}. 

We prove the second part: Let $(Y,S)$ be associated to the $\mathfrak{a}$-adic group $G_{\mathfrak{a}}$, where $\mathfrak{a}=(a_1,a_2,\dots)$. Choose $H$ to be the rational group associated to $\mathfrak{a}$, i.e.\ 
$H = \set{\frac{m}{a_1a_2\cdots a_n}}{m\in\mathbb{Z},\, n\in\mathbb{N}}$. As above we identify $H$ with constant affine functions on $K$ in an obvious way. Let $\mathcal{E}$ be as above, and let $G$ be generated by $H$ and $\mathcal{E}$ 
as above. Let $N$ be any countable torsion-free abelian group, and let $\tilde{G} = G\oplus N$, with $\tilde{G}^+ = G^+\oplus 0$ and order unit $u=1\oplus 0$. It is easily seen that $\mathbb{Q}(\tilde{G},u)\cong \mathbb{Q}(G,1)\cong H$
as ordered groups with distinguished order units. Clearly $\mathrm{Inf}(\tilde{G})  =0\oplus N \cong N$. Also, $S_u(\tilde{G})\cong S_1(G)\cong K$. By Theorem \ref{th:realizing} there exists a Toeplitz flow $(\tilde{X}, \tilde{T})$ such that 
$K^0(\tilde{X}, \tilde{T})\cong \tilde{G}$ (as ordered groups with distinguished order units) and $h(\tilde{T}) = t$. By Proposition \ref{prop:maxfactor} we get that the maximal equicontinuous factor of $(\tilde{X}, \tilde{T})$ is $(Y,S)$,  
and by Theorem \ref{th:realization} we get that $M(\tilde{X}, \tilde{T})\cong K$. Clearly there exists an uncountable family of pairwise non-isomorphic groups $N$ of the above type. The associated simple dimension groups $\tilde{G} = G\oplus N$ are then pairwise non-isomorphic. 
(Note that if $(G_1, G_1^+,u_1) \cong (G_2, G_2^+,u_2)$, then $\mathrm{Inf}G_1 \cong \mathrm{Inf}G_2$.) The corresponding  Toeplitz flows $(\tilde{X},\tilde{T})$ are then pairwise non-strong orbit equivalent by Theorem \ref{th:kgroups}. 
This finishes the proof of Theorem \ref{th:real_simplex}. \ep
\medbreak

\section{Examples}
\label{sec:examples}

We will first state some results that -- combined with our theorems -- will provide a rich source of examples of Toeplitz flows. We first need a definition.

\begin{definition}
 Let $(V,E)$ be a simple Bratteli diagram such that $|V_n|\leq l < \infty$ for all $n$. We say that $(V,E)$ is of finite rank. If $|V_n|= k$ for all $n=1,2,\dots$, we say that $(V,E)$ is of rank $k$.
\end{definition}

\proc{Remark.}
If $(V,E)$ is of finite rank we may telescope $(V,E)$ to get a new Bratteli diagram of rank some $k$. It is easily seen that if $(V,E)$ is of rank $k$, then the dimension group $K_0(V,E)$ has rank $\leq k$.
Furthermore, the state space $S_{[\mathbf{1}]}(K_0(V,E))$ of $K_0(V,E)$ has at most $k$ extreme points, and so is a (finite-dimensional) $m$-simplex for some $m\leq k$.

\begin{theorem}[\cite{DM}]
\label{th:eitheror}
 Let $(V,E,\geq)$ be a properly ordered Bratteli diagram such that $(V,E)$ is of finite rank. Then $(X_{(V,E)},T_{(V,E)})$ is either expansive or it is an odometer.
\end{theorem}

\begin{theorem}
\label{th:zeroentropy}
 Let $(V,E,\geq)$ be a properly ordered Bratteli diagram such that $(V,E)$ has finite rank. Then the entropy of $(X_{(V,E)},T_{(V,E)})$ is zero.
 \end{theorem}

We will prove Theorem \ref{th:zeroentropy} at the end of this section. Combining Theorem \ref{th:eitheror} and Theorem \ref{th:zeroentropy} with Lemma \ref{lem:ers} we get the following result.

\begin{theorem}
\label{th:eitheror2}
 Let $(V,E,\geq)$ be a properly ordered Bratteli diagram with the following two properties:
 \begin{enumerate}[(i)]
  \item $(V,E)$ is of finite rank (and so $K_0(V,E)$ is of finite rank).
  \item $(V,E)$ has the ERS property (and so $\mathbb{Q}(K_0(V,E),[\mathbf{1}])$ is noncyclic).
 \end{enumerate}
Then $(X_{(V,E)},T_{(V,E)})$ is either an odometer or a Toeplitz flow of entropy zero. In particular, if $K_0(V,E)$ is not a (noncyclic) rational group, then $(X_{(V,E)},T_{(V,E)})$ is a Toeplitz flow.
\end{theorem}

There is a partial converse to Theorem \ref{th:eitheror2} due to Handelman \cite[Theorem 8.5]{H}. (Incidentally, the paper \cite{H} treats a more general situation under the assumption that the dimension groups in question have a unique state.) We formulate his result using our terminology and notation. 

\begin{theorem}
\label{th:partconv}
 Let $(G,G^+,u)$ be a simple dimension group with order unit $u$. Assume $\mathbb{Q}(G,u)$ is noncyclic and that $S_u(G)$ is a one-point set. (We say that $G$ has a unique state, the order unit being understood.)
 Assume $\mathrm{rank}(G)=k$. There exists a simple Bratteli diagram $(V,E)$ of rank at most $k+1$ with the ERS property, such that $(G,G^+,u)\cong (K_0(V,E), K_0(V,E)^+,[\mathbf{1}])$ and $(V,E)$ is an ERS realization of $G$ with respect to $\mathbb{Q}(G,u)$. (Cf.\ Remark after Lemma \ref{lem:super}.) Giving $(V,E)$ a proper ordering the associated Bratteli-Vershik
 system is either an odometer or a uniquely ergodic Toeplitz flow of zero entropy.
\end{theorem}

\proc{Remark.}
In particular, the scenario described in Theorem \ref{th:partconv} occurs when $G=H\oplus \mathbb{Z}^m$, where $H$ is a noncyclic rational group, and $G^+=H^+\oplus 0$, with $u=1\oplus 0 $.
Clearly $\mathrm{Inf}(G) = 0\oplus \mathbb{Z}^m \cong \mathbb{Z}^m$, and $\mathbb{Q}(G,u)\cong H$ (as ordered groups with order unit $u$ and 1, respectively).
Clearly $G$ has a unique state and $\mathrm{rank}(G)=m+1$. So by Theorem \ref{th:partconv} there exists an ERS realization $(V,E)$ of $G$ with respect to $H$ such that $(V,E)$ 
has rank at most $m+2$. 

Another example where Theorem \ref{th:partconv} applies is $G=\mathbb{Q}+\mathbb{Q}\alpha\subseteq \mathbb{R}$, where $\alpha$ is an irrational number, and $G$ inherits the ordering from 
$\mathbb{R}$, i.e.\ $G^+ = G\cap \mathbb{R}^+$,
the order unit being 1. Then $\mathbb{Q}(G,1)=\mathbb{Q}$ and $\mathrm{rank}(G)=2$. So by Theorem \ref{th:partconv}, $G$ has an ERS realization $(V,E)$ with respect to $\mathbb{Q}$, with 
$(V,E)$ of rank at most 3. (In fact, in this particular case it suffices with a rank equal to 2, cf.\ Example \ref{ex:indlimit}.) 

We will give explicit examples below illustrating some of these scenarios.
\medbreak

\begin{example}
\label{ex:indlimit}
 Let $G=\mathbb{Q}+\mathbb{Q}\alpha\subseteq\mathbb{R}$ with $G^+$ and order unit 1 as above, and so $\mathbb{Q}(G,1)=\mathbb{Q}$. We may assume without loss of generality that $0<\alpha<\frac{1}{2}$. Let 
  \[\alpha = \frac{1}{a_0} \begin{matrix} \\ +\end{matrix} \frac{1}{a_1}\begin{matrix} \\ +\end{matrix}\frac{1}{a_2}\begin{matrix} \\ +\end{matrix}\cdots\]
 be the continued fraction expansion of $\alpha$. Then $\mathbb{Q}+\mathbb{Q}\alpha$ is order isomorphic to  the inductive limit
 \begin{equation}
  \label{eq:qlim}
  \mathbb{Q}\stackrel{A_1}{\longrightarrow}\mathbb{Q}^2\stackrel{A_2}{\longrightarrow}\mathbb{Q}^2\stackrel{A_3}{\longrightarrow}\mathbb{Q}^2\stackrel{A_4}{\longrightarrow}\cdots \tag{*}
 \end{equation}
where the incidence matrices are 
\[A_1=\begin{bmatrix} a_0 \\ 1 \end{bmatrix},\, A_i = \begin{bmatrix} a_{i-1} & 1 \\ 1 & 0 \end{bmatrix}, \, i\geq 2, \]
and the order unit is the canonical one, i.e.\ represented by $1\in\mathbb{Q}$. This follows from the fact that the dimension group $\mathbb{Z}+\mathbb{Z}\alpha\subseteq \mathbb{R}$ is represented
by \eqref{eq:qlim} with the $\mathbb{Q}$'s replaced by $\mathbb{Z}$'s. (Cf.\ \cite[3.3, Example (ii)]{skau}.) We will assume $a_0 = 1$ and so $A_1=\begin{bmatrix} 1 \\ 1 \end{bmatrix}$. (We make this assumption just for convenience;
nothing essential is changed by this, but the ensuing construction becomes more streamlined.) We will show that the inductive limit in \eqref{eq:qlim} is order isomorphic to the inductive limit
 \begin{equation}
  \label{eq:zlim}
  \mathbb{Z}\stackrel{B_1}{\longrightarrow}\mathbb{Z}^2\stackrel{B_2}{\longrightarrow}\mathbb{Z}^2\stackrel{B_3}{\longrightarrow}\mathbb{Z}^2\stackrel{B_4}{\longrightarrow}\cdots \tag{**}
 \end{equation}
 where the incidence matrices ${B_n}$ have equal row sums -- hence the associated simple Bratteli diagram has the ERS property. The ${B_n}$'s are obtained from the ${A_n}$'s by the following procedure, where we have adapted the 
 construction in the proof of Theorem 11 in \cite{gjerde-johansen} to our setting:
 Let $A_2\begin{bmatrix} 1 \\ 1 \end{bmatrix} = \begin{bmatrix} \alpha_1 \\ \alpha_2 \end{bmatrix}$. Then $J_2'A_2\begin{bmatrix} 1 \\ 1 \end{bmatrix} = \begin{bmatrix} 1 \\ 1\end{bmatrix}$, where $J_2'$ is the diagonal $2\times 2$ matrix 
 $\mathrm{diag}(\alpha_1^{-1}, \alpha_2^{-1})$. Let $m_2$ be the least common multiple of the denominators of the entries of $J_2'A_2$, and let $k_2=2m_2$. Set $J_2 = k_2J_2'$. Then all entries of $J_2A_2$ are in $\mathbb{N}$ and are 
 divisible by $2$. Furthermore, $J_2A_2\begin{bmatrix} 1 \\ 1 \end{bmatrix} = \begin{bmatrix} k_2 \\ k_2\end{bmatrix}$. Assume we have constructed $J_2, \dots,J_{n-1}$. Let $J_n'$ be the appropriate diagonal matrix over $\mathbb{Q}^{+}$ such that 
 $J_n'A_nJ_{n-1}^{-1}\begin{bmatrix} 1 \\ 1 \end{bmatrix} = \begin{bmatrix} 1 \\ 1\end{bmatrix}$. Let $m_n$ be the least common multiple of the denominators of the entries of $J_n'A_nJ_{n-1}^{-1}$, and let $k_n=nm_n$. Set $J_n = k_nJ_n'$.
 Then all entries of $J_nA_nJ_{n-1}^{-1}$ are in $\mathbb{N}$ and are divisible by $n$. Furthermore, $J_nA_nJ_{n-1}^{-1}\begin{bmatrix} 1 \\ 1 \end{bmatrix} = \begin{bmatrix} k_n \\ k_n\end{bmatrix}$, where we observe that $n$ is 
 a divisor of $k_n$. Set $B_{n} = J_nA_nJ_{n-1}^{-1}$. Setting $J_0=\mathrm{id}$, $J_1 = \mathrm{id}$, we get the commutative diagram
\begin{center}
 \includegraphics{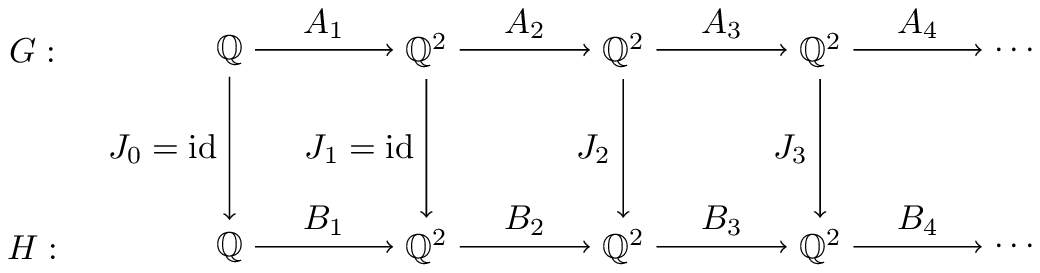}
\end{center}
which establishes an order isomorphism between the two associated inductive limits $G$ and $H$ (respecting the distinguished order units). Now $H$ is order isomophic to the dimension group associated to \eqref{eq:zlim} since the latter is a divisible group (because $n$ is a divisor of each entry of $B_n$). We conclude that $G= \mathbb{Q}+\mathbb{Q}\alpha\cong K_0(V,E)$ (as ordered groups with distinguished order units), where $(V,E)$ is the Bratteli diagram associated 
to \eqref{eq:zlim}. Since $B_n\begin{bmatrix}1\\1\end{bmatrix} = \begin{bmatrix}k_n\\k_n\end{bmatrix}$, we see that $(V,E)$ has the ERS property. Furthermore, since $n$ divides $k_n$ we observe that $(V,E)$
is a rank 2 ERS realization of $G$ with respect to $\mathbb{Q}(G,1)\,(=\mathbb{Q})$. 

Finally, any proper ordering of $(V,E)$ will yield a Bratteli-Vershik (BV-) system which is a uniquely ergodic Toeplitz system of entropy zero. An intriguing question is if by considering all proper orderings there arises an uncountable family 
of pairwise non-conjugate BV-systems? These BV-systems are of course all strong orbit equivalent, having the same $K^0$-group $\mathbb{Q}+\mathbb{Q}\alpha$.
\end{example}

\begin{example}
\label{ex:2sym}
\emph{2-symmetric Bratteli diagrams $(V,E)$ and their associated simple dimension groups.}
\begin{figure}
\includegraphics{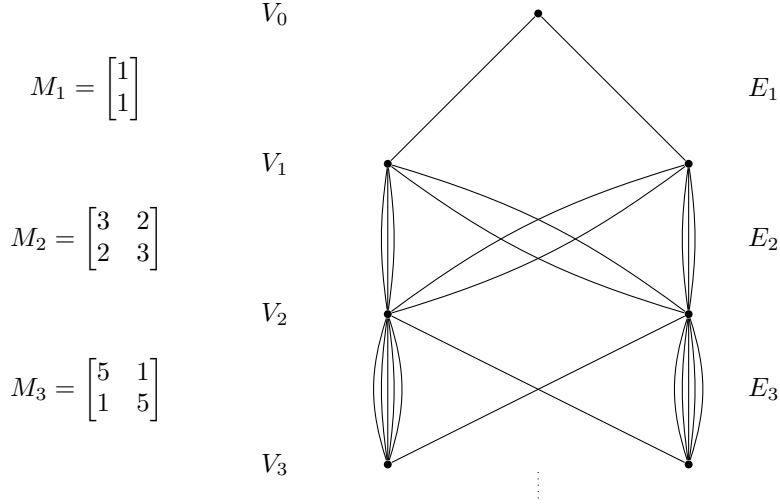}
\caption{A 2-symmetric Bratteli diagram}
\label{fig:2sym}
\end{figure}
Let $|V_n|=2$ for all $n\geq 1$ and let the incidence matrix $M_n$ between levels $n-1$ and $n$ be 2-symmetric, i.e.\ of the form 
\[M_n = \begin{bmatrix} l_n & k_n \\ k_n & l_n \end{bmatrix}, \, \text{ where } 1\leq k_n < l_n\]
for all $n\geq 2$, and $M_1 = \begin{bmatrix} 1\\1\end{bmatrix}$. (See Figure \ref{fig:2sym}.) So $\mathrm{rank}(V,E)=2$ and $(V,E)$ has the ERS-property (as well as the ECS property, i.e.\ equal column sums). (The 2-symmetric Bratteli diagrams were studied in \cite{FM} in 
connection with classifying the symmetries of UHF C$^\ast$-algebras. Cf.\ also \cite[Chapter III, Section 7.7.4]{B}.) It will be convenient to write for $n\geq 2$:
\[l_n=\frac{q_n+r_n}{2}, \quad k_n = \frac{q_n-r_n}{2},\, \text{ where } q_n \text{ and } r_n \text{ have the same parity,}\]
and so $q_n=l_n+k_n$, $r_n=l_n-k_n$. By a simple computation we get that 
\[M_nM_{n-1}\cdots M_2 = \frac{1}{2}\begin{bmatrix}s_n+t_n&s_n-t_n\\s_n-t_n&s_n+t_n\end{bmatrix},\, \text{ where } s_n=\prod_{i=2}^n q_i,\, t_n = \prod_{i=2}^n r_i.\]
We note that $s_n$ and $t_n$ have the same parity and that 
\[\frac{s_n}{t_n} = \prod_{i=2}^n \frac{q_i}{r_i} \nearrow \prod_{i=2}^\infty \frac{q_i}{r_i} = \alpha.\]
Two scenarios can occur. Firstly, if $\alpha< \infty$ then the state space $S_{[{\bf 1}]}(K_0(V,E))$ of $K_0(V,E)$ has two extreme points. Specifically,
\[K_0(V,E) = \left\lbrace \left.\left(\frac{k}{2s_n},\frac{l}{2t_n}\right) \right| k,l\in\mathbb{Z} \text{ with the same parity, } n\in\mathbb{N} \right \rbrace\subseteq \mathbb{Q}^2\subseteq \mathbb{R}^2.\]
$K_0(V,E)^+$ consists of those elements in $K_0(V,E)$ that lie inside the cone determined by the two half-lines with slopes $\alpha$ and $-\alpha$ respectively. The canonical order unit $[{\bf 1}]$ of $K_0(V,E)$ is represented by $(1,0)\, \left(=\left(\frac{2}{2},0\right)\right)$ and the two extreme states are determined by projecting orthogonally to the two lines with slopes $\alpha$ and $-\alpha$, respectively. We note that $K_0(V,E)$ does not split as a direct sum of two (noncyclic) 
rational groups, although it "almost" does. (More on that when describing the unique state case below.) Furthermore, $\mathbb{Q}(K_0(V,E),[{\bf 1}])$ is the subgroup of $K_0(V,E)$ represented by
\[\left\lbrace \left.\left(\frac{k}{2s_n},0\right) \right| k\in\mathbb{Z} \text{ even, } n\in\mathbb{N} \right \rbrace\]
and so is order isomophic to the noncyclic rational group $G(m)$, where $m=\prod_{i=1}^\infty q_i$. 

The second scenario occurs when $\alpha=\infty$. Then $S_{[{\bf 1}]}(K_0(V,E))=\lbrace \tau\rbrace$ is a one-point set. The unique state $\tau$ is intimately related to the Perron-Frobenius (P-F) eigenvalues and eigenvectors of the incidence matrices 
$\lbrace M_n\rbrace_{n=2}^\infty$. In fact, $q_n$ is the P-F eigenvalue og $M_n$ with left (right) eigenvector $(1,1)$ ($(1,1)^\mathrm{tr}$). The other eigenvalue
is $r_n$ with left (right) eigenvector $(1,-1)$ ($(1,-1)^\mathrm{tr}$). The state $\tau$ is determined by 
\[\tau(v_n)=\frac{1}{q_1q_2\cdots q_n}(1,1)\stackrel{\rightarrow}{v_n}\]
where $v_n$ denotes the element in $K_0(V,E)$ represented by the column vector $\stackrel{\rightarrow}{v_n}\in \mathbb{Z}^2$ at level $n$ of the Bratteli diagram, and where we for normalization purposes have set $q_1=2$ so that $\tau([{\bf 1}])=1$. In particular we get that 
\[\tau(K_0(V,E)) = \left \lbrace \left. \frac{m}{q_1q_2\cdots q_n}\right| m\in\mathbb{Z},n\in\mathbb{N}\right\rbrace = Q \, (\subseteq \mathbb{Q}).\]
Furthermore, $\tau(H)= \left \lbrace \left. \frac{m}{q_1q_2\cdots q_n}\right| m\in\mathbb{Z},n\in\mathbb{N}\right\rbrace$ where $H=\mathbb{Q}(K_0(V,E),[{\bf 1}])$, since $H$ is represented as "constant" vectors at each level of $(V,E)$, cf.\ Lemma \ref{lem:super}. We notice that $\tau(H)=2\tau(K_0(V,E))=2Q$. Now we focus on the interesting special case where $r_n=1$, for all $n\geq 2$. Hence $l_n=k_n+1$, and $q_n$ is odd for all $n\geq 2$. In particular, $Q=\tau(K_0(V,E))$ is not 2-divisible. Now $\mathrm{ker}(\tau)=\mathrm{Inf}(K_0(V,E))$ and $K_0(V,E)^+=\tau^{-1}(\mathbb{R}^+-\{0\})$. We observe that $\mathrm{ker}(\tau)\cong \mathbb{Z}$ since $\tau(v_n)=0$ for all $n\geq 1$, where $\stackrel{\rightarrow}{v_n}=(1,-1)^\mathrm{tr}$ (cf.\ the notation above), and $M_n(1,-1)^\mathrm{tr} = (1,-1)^\mathrm{tr}$ for all $n\geq 2$. Thus we get the short exact sequence 
\[0\longrightarrow \mathbb{Z}\longrightarrow K_0(V,E) \stackrel{\tau}{\longrightarrow} Q \longrightarrow 0\]
However, there is no splitting $K_0(V,E)=Q\oplus\mathbb{Z}$, with $\tau$ the projection onto $Q$, the reason being that $Q$ is not 2-divisible. In fact, if such a splitting occurred then clearly $Q$ would be equal to the rational subgroup $H$. This would imply that $\tau(H)=Q$ and, combining this with $\tau(H)=2Q$, we would get that $Q$ is 2-divisible which is a contradiction. (For more details, cf. \cite[Section 8]{H}.)
\end{example}

\proc{Remark.}
By our results we know that all proper orderings of a 2-symmetric Bratteli diagram yield Bratteli-Vershik maps that are Toeplitz flows. One can show that when $r_n=1$ ($n\geq $) and $(V,E)$ is given the left-right ordering (meaning that we order the edges ranging at the same vertex from left to right), then the Toeplitz sequence that corresponds to the (unique) minimal path is regular (cf. Definition \ref{def:regular}).
\medbreak

\proc{Proof of Theorem \ref{th:zeroentropy}}
By telescoping we may assume that $k = |V_1|=|V_2|= \cdots (=\mathrm{rank}(V,E))$, where $V = \{V_0,V_1,V_2,\dots\}$, $E = \{E_1, E_2, \dots\}$. By Theorem \ref{th:eitheror} we may assume that
$(X_{(V,E)}, T_{(V,E)})$ is expansive, since if not it would be an odometer which has zero entropy. This implies in particular that we may assume that $k\geq 2$. By telescoping the initial part of $(V,E)$ we may assume that $(X_{(V,E)}, T_{(V,E)})$ is conjugate
to $(X_1,S_1)$ (and hence conjugate to $(X_n,S_n)$ for all $n\geq 1$) by the map $\deff{\pi_1}{X_{(V,E)}}{X_1\,(=\pi_1(X_{(V,E)})\subseteq E_1^\mathbb{Z}\,(=P_1^\mathbb{Z}))}$ defined by 
$\pi_1(x)=\left(\tau_1(T_{(V,E)}^nx)\right)_{n=-\infty}^\infty$, where $\deff{\tau_1}{X_{(V,E)}}{E_1\,(=P_1)}$ is the truncation map. (Cf.\ Proposition \ref{prop:expansive-conj} and the description
of key constructions associated to $(V,E,\geq)$, as well as notation and terminology preceding the proposition.) For $v\in V_n$ we consider all the paths, say $\{p_1, p_2, \dots, p_l\}$ from 
$v_0\in V_0$ to $v$, which will be a subset of $P_n$. Let those paths have the ordering $p_1<p_2<\cdots<p_l$, say, in the induced lexicographic ordering. Then $w(v) = \tau_1(p_1)\tau_1(p_2)\cdots\tau_1(p_l)$
will be a word over $E_1$. Let $W_n=\{w(v)\,|\, v\in V_n\}$. Notice that $|W_n|=k$ for all $n$. Define $(X_{W_n},S_{W_n})$ to be the subshift of the full shift on $E_1^\mathbb{Z}$, where $X_{W_n}$ is the set of all 
bisequences formed by concatenation of words in $W_n$. Now we observe that $E_1^\mathbb{Z}\supseteq X_{W_1}\supseteq X_{W_2}\supseteq \cdots \supseteq X_1$. In fact, this is an immediate consequence of how the Vershik
map is defined. 

Recall that the entropy of a subshift $(X,T)$ is $h(T)=\lim_{q\rightarrow \infty} \frac{1}{q}\log{|B_q(X)|}$, where $B_q(X)$ is the set of words of length $q$ occurring in $X$. (Cf. \cite[Theorem 7.13]{walters}.)
Clearly $h(S_1)\leq h(S_{W_n})$ for all $n$, and so it suffices to show that $h(S_{W_n})\rightarrow 0$ as $n\rightarrow \infty$. Let $l_n$ be the length of the shortest word in $W_n$. (Clearly $l_n\rightarrow\infty$ as $n\rightarrow \infty$.)
Assume $w\in X_{W_n}$ is a subword of length $m$ of a concatenation of $s$ words from $W_n$. Then we easily see that 
\[s\leq \ceil{\frac{m}{l_n}}+1\]
where $\ceil{x}$ is the smallest integer which is larger or equal to $x$. There is at most $k^s$ different ways to concatenate $s$
words in $W_n$, and so we get 
\[|B_m(W_n)| \leq \sum_{s=0}^{\ceil{\frac{m}{l_n}}+1}k^s = \frac{k^{\ceil{\frac{m}{l_n}}+2}-1}{k-1}\leq k^{\ceil{\frac{m}{l_n}}+2} \leq k^{\frac{m}{l_n}+3}.\]
This implies that $\frac{1}{n}\log{|B_m(W_n)|}\leq \frac{\frac{m}{l_n}+3}{n}\log{k}\rightarrow 0$ as $n \rightarrow \infty$. Hence 
$h(S_{W_n})\rightarrow 0 $ as $n \rightarrow \infty$, thus completing the proof.
\ep
\medbreak

\section*{Acknowledgement}
I would like to thank David Handelman and Fumiaki Sugisaki for helpful private communications. Furthermore, some of the results in this paper depend upon earlier results obtained by them. Finally, I would like to extend my sincere gratitude to Christian Skau for suggesting the basic problems considered in this paper and for his constant encouragement and valuable advice.

\addcontentsline{toc}{section}{Bibliography}
\bibliographystyle{amsalpha}

\end{document}